\pdfoutput=1
\RequirePackage{ifpdf}
\ifpdf 
\documentclass[pdftex]{sigma}
\else
\documentclass{sigma}
\fi

\makeatletter
\newcommand*\bigcdot{\mathpalette\bigcdot@{.5}}
\newcommand*\bigcdot@[2]{\mathbin{\vcenter{\hbox{\scalebox{#2}{$\m@th#1\bullet$}}}}}
\makeatother

\def\R{{\mathbb R}}

\newtheorem{thm}{Theorem}[section]

\newtheorem{prop}[thm]{Proposition}
\newtheorem{cor}[thm]{Corollary}
{ \theoremstyle{definition}
\newtheorem{rem}[thm]{Remark}
\newtheorem{defn}[thm]{Definition}}

\def\fracz#1#2{#1/#2}

\def\ro#1{{\rm #1}}

\def\Bbb#1{{\mathbb#1}}

\def\R{\Bbb R}

\def\Rx{\R\mkern1mu^}

\def\Rn{\Rx n}

\def\AG#1{\ro{A}(#1)} 
\def\AGR#1{\AG{#1,\R}} 
\def\SA#1{\ro{SA}(#1)} 
\def\SAR#1{\SA{#1,\R}} 
\def\GL#1{\ro{GL}(#1)} 
\def\GLR#1{\GL{#1,\R}} 
\def\SL#1{\ro{SL}(#1)} 
\def\SLR#1{\SL{#1,\R}}  

\def\semidirect{\ltimes}

\def\ip#1#2{\langle \,{#1}\,,{#2}\,\rangle}
\def\nnorm#1{ \|\,{#1}\,\|}

\numberwithin{equation}{section}

\begin{document}
\allowdisplaybreaks

\newcommand{\arXivNumber}{2003.13842}

\renewcommand{\PaperNumber}{093}

\FirstPageHeading

\ShortArticleName{Feature Matching and Heat Flow in Centro-Affine Geometry}

\ArticleName{Feature Matching and Heat Flow\\ in Centro-Affine Geometry}

\Author{Peter J. OLVER~$^\dag$, Changzheng QU~$^\ddag$ and Yun YANG~$^\S$}

\AuthorNameForHeading{P.J.~Olver, C.~Qu and Y.~Yang}

\Address{$^\dag$~School of Mathematics, University of Minnesota, Minneapolis, MN 55455, USA}
\EmailD{\href{mailto:olver@umn.edu}{olver@umn.edu}}
\URLaddressD{\url{http://www.math.umn.edu/~olver/}}

\Address{$^\ddag$~School of Mathematics and Statistics, Ningbo University, Ningbo 315211, P.R.~China}
\EmailD{\href{mailto:quchangzheng@nbu.edu.cn}{quchangzheng@nbu.edu.cn}}

\Address{$^\S$~Department of Mathematics, Northeastern University, Shenyang, 110819, P.R.~China}
\EmailD{\href{mailto:yangyun@mail.neu.edu.cn}{yangyun@mail.neu.edu.cn}}

\ArticleDates{Received April 02, 2020, in final form September 14, 2020; Published online September 29, 2020}

\Abstract{In this paper, we study the differential invariants and the invariant heat flow in centro-affine geometry, proving that the latter is equivalent to the inviscid Burgers' equation. Furthermore, we apply the centro-affine invariants to develop an invariant algorithm to match features of objects appearing in images. We show that the resulting algorithm compares favorably with the widely applied scale-invariant feature transform (SIFT), speeded up robust features (SURF), and affine-SIFT (ASIFT) methods.}


\Classification{53A15; 53A55}

\section{Introduction}

The main objective in this paper is to study differential invariants and invariant curve flows~-- in particular the heat flow~-- in centro-affine geometry. In addition, we will present some basic applications to feature matching in camera images of three-dimensional objects, comparing our method with other popular algorithms.

Affine differential geometry is based on the Lie group $\AGR n = \GLR n \semidirect \Rn$ consisting of affine transformations $x \longmapsto Ax+b$, $A\in \GLR n$, $b\in \Rn$ acting on $x \in \Rn$. Basic references include Nomizu and Sasaki~\cite{ns} and Simon~\cite{sim}. Keep in mind that, in most of the literature, the term ``affine geometry'' usually refers to ``equi-affine geometry'', in which one restricts to the subgroup $\SAR n = \SLR n \semidirect \Rn$ of volume-preserving affine transformations. A key issue is to study the resulting invariants associated with submanifolds $M \subset \Rn$. In particular, the classical theory for equi-affine hypersurfaces was developed by Blaschke and his collaborators,~\cite{Blaschke}.

\emph{Centro-affine differential geometry} refers to the geometry induced by the general linear group $x \longmapsto Ax$, $A\in \GLR n$, $x \in \Rn$, which is the subgroup of the affine transformation group that keeps the origin fixed. Similarly, \emph{centro-equi-affine differential geometry} refers to the subgroup $\SLR n$ of volume-preserving linear transformations. These cases are usually discussed in passing in books that are devoted to (equi-)affine geometry~\cite{ns}.

Several methods have been developed to construct differential invariants and other invariant quantities in such Klein geometries, \cite{gw, olv-5,olv-6,pk,wil}. In particular, invariants can be straightforwardly and algorithmically obtained by the method of equivariant moving frames introduced in~\cite{fo-2}. They play a prominent role in the study of the geometric properties, including equivalence and symmetry, of curves, surfaces, and more general submanifolds, as well as invariant geometric flows,~\cite{olv-41}, with many applications to computer vision and image processing.

The term ``invariant submanifold flow'' refers to the motion of a curve or surface by a prescribed partial differential equation that admits an underlying transformation group as a symmetry group, e.g.~the Euclidean group of rigid motions (translations and rotations). Invariant curve and surface flows arise in an impressive range of applications, including geometric optics, elastodynamics, computer vision, visual tracking and control, interface dynamics, thermal grooving, and elsewhere. A celebrated example is the curve shortening flow (CSF), in which a plane curve moves in its normal direction in proportion to its curvature. The CSF was first introduced by Mullins \cite{mul} as a model for the motion of grain boundaries. It was shown by Gage and Hamilton~\cite{gh} and Grayson \cite{gra-1,gra-2} that a simple closed Jordan curve will remain simple when evolving under the CSF, first becoming convex and then shrinking to a point in finite time while becoming asymptotically circular, often referred to as a ``circular point''. These results were an important preliminary step on the road to Hamilton's celebrated analysis of the higher dimension counterparts such as mean curvature flow and Ricci flow; the latter was extended by Perelman to in his famous solution to the Poincar\'e conjecture,~\cite{ds}. The corresponding affine curve shortening flow (ACSF) was introduced and studied by Angenent, Sapiro, and Tannenbaum \cite{ast, st-2}, and it was proved that a closed convex curve when evolves according to the ACSF will shrink to an ``elliptical point''. See also \cite{and,cz} and the references therein for further developments. In computer vision, \cite{olv-6, ost-1,ost-2}, Euclidean curve shortening and its equi-affine counterpart have been successfully applied to image denoising and segmentation and are actively used in practical computer implementations, both academic and commercial. Euclidean-invariant three-dimensional curve flows include the integrable vortex filament flow appearing in three-dimensional fluid dynamics, \cite{Hasimoto,LangPer}, while mean curvature and Willmore flows of surfaces have been the subject of extensive analysis and applications, \cite{dde,KuSc}. More recently, similar results were obtained for the heat flow in centro-equi-affine geometry \cite{wwq}. Heat flows in more general Klein geometries were proposed~\cite{olv-41,ost-1}.

In this paper, we are interested in the heat flow in centro-affine geometry. Interestingly, we find that the heat flow for the centro-affine curvature $\kappa(t,s)$ is equivalent to the well-known \emph{inviscid Burgers' equation} $\kappa_t=\kappa\kappa_s$. This result is in contrast to the behavior of heat flows in Euclidean, equi-affine, and centro-equi-affine geometries, which yield second order nonlinear parabolic equations for the associated curvature invariant.

A challenging problem arising in computer vision and pattern recognition, is feature matching under viewpoint changes between different images. Image and feature matching has wide range of applications, including robotic vision, medical image registration, 3D reconstruction, optical character recognition, object classification, content-based image retrieval, and so on. Traditional methods, such as scale invariant feature transform (SIFT)~\cite{low} and speeded up robust features (SURF)~\cite{btv}, have excellent performance and high precision.
 However, when the images have less texture complexity and color diversity, it is not easy to extract and describe the feature points.
 Another drawback to these detectors is that they are only invariant under the planar Euclidean group consisting of rigid motions (rotation, translation, and, possibly, reflection).
 In fact, the apparent deformations of three-dimensional objects caused by changes of the camera
 position can be locally approximated by affine maps, and hence, during the image matching process, the Euclidean transformation
 group should be extended to the affine transformation group by including stretching and skewing transformations.
 Applications of equi-affine and affine invariance to image processing can be found, for instance, in \cite{ms,mts,tot}. An affine-invariant extension of the SIFT algorithm (ASIFT) has been proposed
 in \cite{rdm,ym}, which detects feature points in two images that are so related by
 simulating many affine transformations of each image and performing the SIFT algorithm between all
 image pairs. Going beyond affine-invariant detectors, projective invariance and moving frame-based signatures have also been successfully applied in computer vision applications~\cite{hh}. Recently, Damelin, Ragozin and Werman investigated the best uniform approximation to a continuous function under affine transformations, which has applications in the rapid rendering of computer graphics~\cite{drw}.

 In contrast to the above mentioned methods, a centro-affine invariant detection method offers the following features:
 \begin{enumerate}\itemsep=0pt
\item[(1)] It depends on the centro-affine differential invariants of smooth closed curves, in contrast to the ordinary discrete mode matching methods.
\item[(2)] Centro-affine differential invariants involve lower order derivatives of the curve para\-met\-ri\-za\-tion, and hence are more accurate and less error prone than their fully affine counterparts.
\item[(3)] In some situations, centro-affine invariance is equivalent to fully affine invariance if we can find a pair of exact corresponding points (a point-correspondence) with respect to an affine transformation.
 More precisely, for every closed curve, its barycenter can serve as that point-correspondence (or local origin) for its local centro-affine invariants.
\item[(4)] The method relies solely on edge detection, and hence can be applied to untextured images.
\end{enumerate}

 The remainder of this paper is organized as follows. In Section~\ref{sec-back}, we provide a brief review discussion on the moving frame method and differential invariants. In Section~\ref{sec-cai}, a classification for the planar curves with constant centro-affine curvatures is provided. In Section~\ref{section4}, we study the centro-affine invariant heat flow. In Section~\ref{section5}, an application of the centro-affine invariants in the matching of images obtained by cameras is discussed. Finally, Section~\ref{section6} contains some concluding remarks on this work.

\section{Preliminaries}\label{sec-back}

 \subsection{Moving frame}\label{mvf}

 Let us first review basic facts on the method of equivariant moving frames introduced by Fels and the first author~\cite{fo-2, olv-6}.
 Assume $G$ is an $r$-dimensional Lie group acting smoothly on an $m$-dimensional manifold $M$:
 \begin{equation*}
 G\times M \rightarrow M,\qquad h\cdot(g\cdot z)=(hg)\cdot z.
 \end{equation*}
 A {\it right equivariant moving frame} is defined as a smooth map $\rho\colon M\rightarrow G$,
 that is equivariant with respect to the action on $M$ and the inverse right action of $G$ on itself;
 explicitly,
 \begin{equation}\label{rmf}
 \rho\colon \ M\rightarrow G \quad\hbox{satisfies} \quad \rho(g\cdot z)=\rho(z)g^{-1}.
\end{equation}
The existence of a (local) moving frame requires that the group act freely and regularly on $M$. The regularity is a global condition and does not play a role in any of the applications to date. In many cases, one gets by with a locally free action, in which case the resulting moving frame is locally equivariant, meaning that \eqref{rmf} holds for group elements $g$ sufficiently close to the identity.

 Given local coordinates $z=(z_1,\ldots,z_m)$ on $M$,
 let $w(g,z)=g\cdot z$ be the formulae for the transformed coordinates under the group transformation.
 The right moving frame $g=\rho(z)$ associated to the coordinate cross-section $K=\{z_1=c_1,\ldots,z_r=c_r\}$ is found by solving the {\it normalization equations}
 \begin{equation}\label{mo-equa}
 w_1(g,z)=c_1,\ \ldots, \ w_r(g,z)=c_r,
 \end{equation}
 for the group parameters $g = (g_1,\ldots,g_r)$ in terms of the coordinates $z=(z_1,\ldots,z_m)$. The equivariant
 moving frames can be obtained by choosing a cross-section $K$ and solving for the group element $g = \rho(z)$ that
 takes $z$ to a point $\rho(z) \cdot z \in K$, known as the ``canonical form'' of~$z$.
 The coordinates of the canonical form provide a complete system of non-constant invariants.

 \begin{thm}
 If $g=\rho(z)$ is the moving frame obtained by solving equation~\eqref{mo-equa}, then
 \begin{equation*}
 I_1(z)=w_{r+1}(\rho(z)\cdot z), \ \ldots , \ I_{z-r}(z)=w_m(\rho(z)\cdot z),
 \end{equation*}
 form a complete system of functionally independent invariants.
 \end{thm}

 \begin{defn} The invariantization of a scalar function $F\colon M\rightarrow \R$
 with respect to a right moving frame is the invariant function $I=\iota(F)$ defined by $I(z)=F(\rho(z)\cdot z)$. In particular, if~$I$ is any invariant function, then $I=\iota(I)$.
 \end{defn}

 Thus, invariantization defines a canonical projection from the algebra of (smooth) functions to the algebra of invariant functions that respects all algebraic operations.

 \begin{defn}Given a smooth manifold $M$ of dimension $m$ and an integer $1 \leq p<m$,
 the $k$-th order jet bundle $J^k=J^k(M,p)$ is a fiber bundle over $M$,
 such that the fiber of a point $z\in M$ consists of the set of equivalence classes of $p$-dimensional submanifolds of $M$ under the equivalence relation of $k$-th order contact at the point $z$.
 \end{defn}

The \emph{regular subset} of the jet bundle is where the action is (locally) free and regular, which is non-empty when the order~$k$ is sufficiently large. A~(local) moving frame of order $k$ can then be constructed on the regular subset. A $p$-dimensional submanifold is called \emph{regular} at order $k$ if its jet belongs to the regular subset, and hence is in the domain of the moving frame map. See~\cite{olv-sinmf} for an algebraic characterization of totally singular submanifolds, meaning those whose jets are singular at all orders. In centro-affine geometry, the straight lines are totally singular.

 Assume the manifold $M$ has local coordinates $z = \big(x^1,\ldots,x^p,u^1,\ldots,u^q\big)$ in some neighborhood
 where the regular submanifold $S$ can be represented as a graph $u=u(x)$.
 The {\it fundamental differential invariants} are obtained by invariantization of the individual jet coordinate functions,
 \begin{equation*}
 H^i=\iota\big(x^i\big),\qquad I^\alpha_J=\iota\big(u^\alpha_J\big),\qquad i=1,\ldots,p, \quad \alpha=1,\ldots,q, \quad \#J\geq 0.
 \end{equation*}

One can further apply the invariantization process to differential forms by the same procedure. First transform the differential form by acting on it by a general group element $g$ and then invariantization by replacing the group parameters by their expressions $g = \rho(z)$ in terms of the moving frame. Details can be found in~\cite{fo-2}. In particular the (horizontal components of) the invariantized basis horizontal one-forms ${\rm d}x^i$ give the fundamental (contact-)invariant one-forms, which we denote by\footnote{We ignore all contact forms,~\cite{olv-eis}, which is the meaning of the $H$ subscript.}
 \begin{equation}\label{omega}
\omega ^i = \iota\big({\rm d}x^i\big)_H,\qquad i = 1,\ldots,p.
\end{equation}
In curve geometries, $\omega = \iota({\rm d}x)_H$ is the invariant arc length element, usually denoted as ${\rm d}s$. The corresponding dual \emph{invariant differential operators} are denoted by $\mathcal D_1,\ldots, \mathcal D_p$, and can be directly obtained by substituting the moving frame formulas for the group parameters into the corresponding implicit differentiation operators used to produce the prolonged group actions. The invariant differential operators map differential invariants to differential invariants, and hence can be iteratively applied to generate the higher order differential invariants.

\begin{thm} If the moving frame has order $n$, then the set of fundamental differential invariants
\begin{equation*}
 \mathcal{I}^{(n+1)} =
 \big\{H_i,I^\alpha_J\, |\, i=1,\ldots p,\, \alpha=1,\ldots,q,\, \#J\leq n+1\big\}
 \end{equation*}
of order $\leq n+1$ forms a generating set, meaning that all other differential invariants can be obtained by invariant differentiation.
 \end{thm}

In many cases, $\mathcal{I}^{(n+1)}$ does not form a minimal generating set, owing to the existence of syzygies (algebraic relations) among the differentiated invariants. Nevertheless, these syzygies and indeed the entire structure of the differential invariant algebra can be completely determined by the moving frame calculus, using the powerful recurrence formulae. See \cite{fo-2,olv-6} for details.

\subsection{Centro-affine differential invariants for plane curves}

Let us now implement the moving frame calculation, based on Section \ref{mvf}, for the centro-affine group $\AGR2$ acting on plane curves.
 In this case the group is acting on $J^5\big(\R^2, 1\big)$. We represent the planar centro-affine action explicitly in a slightly more convenient form:
 \begin{equation*}
 \left(
 \begin{matrix}
 u \\
 v
 \end{matrix}
 \right)
 =
 \lambda\left(
 \begin{matrix}
 \alpha & \beta \\
 \gamma & \delta
 \end{matrix}
 \right)
 \left(
 \begin{matrix}
 x \\
 y
 \end{matrix}
 \right),
 \qquad\hbox{where} \quad
 \det\left(
 \begin{matrix}
 \alpha & \beta \\
 \gamma & \delta \\
 \end{matrix}
 \right)=1,
 \quad
 \lambda\neq 0.
 \end{equation*}
 By a direct computation,
 the prolonged centro-affine transformations up to order $4$ are given by
 \begin{gather*}
 u =\lambda\alpha x+\lambda\beta y, \qquad v=\lambda\gamma x+\lambda\delta y,\\
 v_u =\frac{\gamma+\delta y_x}{\alpha+\beta y_x}, \qquad v_{uu}=\frac{y_{xx}}{\lambda(\alpha+\beta y_x)^3},\\
 v_{uuu} =\frac{(\alpha+\beta y_x)y_{xxx}-3\beta y^2_{xx}}{\lambda^2(\alpha+\beta y_x)^5},\\
 v_{uuuu} =\frac{(\alpha+\beta y_x)^2y_{xxxx}-10\beta(\alpha+\beta y_x)y_{xx}y_{xxx}+15\beta^2y^3_{xx}}{\lambda^3(\alpha+\beta y_x)^7}.
 \end{gather*}
 Let $\epsilon=\operatorname{sign} [(y_{xx}/(xy_x-y) ]$. Further, after possibly reparametrizing or applying a centro-affine transformation, we can specify $\operatorname{sign}(y-xy_x)=1$. To construct a moving frame, we use the cross-section normalization
 \begin{equation*}
 u=0,\qquad v=1,\qquad v_u=0,\qquad v_{uu}=-\epsilon.
 \end{equation*}
 Solving for the group parameters yields
 \begin{gather}
 \lambda^4=\frac{\epsilon y_{xx}}{(xy_x-y)^3},\qquad \alpha=\lambda y,\qquad \beta=-\lambda x,\nonumber\\ \gamma=-\frac{y_x}{\lambda(y-xy_x)},\qquad \delta=\frac{1}{\lambda(y-xy_x)},\label{mvA2}
 \end{gather}
which prescribes the right-equivariant moving frame. Invariantizing the horizontal one-form
\begin{equation*}
{\rm d}u_H = (\lambda\alpha +\lambda\beta y_x) {\rm d}x
\end{equation*}
 by substituting the moving frame formulae \eqref{mvA2} produces the centro-affine arc-length element
 \begin{equation}\label{aff_arcl}
 {\rm d}s=\iota({\rm d}x)_H = {\operatorname{sign}}(y-xy_x)\sqrt{\epsilon \frac{y_{xx}}{xy_x-y}}\ {\rm d}x=\sqrt{\epsilon \frac{y_{xx}}{xy_x-y}} {\rm d}x,
 \end{equation}
 By a direct calculation, we produce the fundamental differential invariant
 \begin{gather}
 \kappa = \iota(y_{xxx})= {\operatorname{sign}}(y-xy_x)\frac{3xy^2_{xx}+(y-xy_x)y_{xxx}}{y^2_{xx}}\sqrt{\epsilon\,\frac{y_{xx}}{xy_x-y}}\nonumber\\
 \hphantom{\kappa = \iota(y_{xxx})}{} = \frac{3xy^2_{xx}+(y-xy_x)y_{xxx}}{y^2_{xx}}\sqrt{\epsilon\,\frac{y_{xx}}{xy_x-y}}, \label{aff_curv}
 \end{gather}
 which we identify as the \emph{centro-affine curvature}.
 Higher order invariants are all obtained by invariant differentiation of~$\kappa $ with respect to the centro-affine arc length \eqref{aff_arcl}, and so a complete list of differential invariants is given by $\kappa , \kappa _s, \kappa _{ss}, \ldots$. In particular, one can show, either by direct calculation or by using the recurrence formulae~\cite{fo-2}, that
 \[ \iota(y_{xxxx})=\kappa_s+\frac{3}{2}\kappa^2-3,\qquad\iota(y_{xxxxx})=\kappa_{ss}+5\kappa\kappa_s+3\kappa^3-16\kappa,\]
 and so on.

It is also of interest to obtain the formulas for the centro-affine curvature and arc length for a general parametrized curve.
Consider a smooth curve parameterized by
 \begin{equation*}
 \mathbf{x}(p)= (x(p), y(p) )^{\mathrm{T}},
 \end{equation*}
 where $x(p)$, $y(p)$ are smooth functions of the parameter $p$ defined over a certain interval $I$,
 and the superscript ``${}^\mathrm{T}$'' represents the transpose of a vector or matrix. We use dots to denote derivatives with respect to the parameter $p$. In particular $\dot{\mathbf{x}}=\fracz{{\rm d}\mathbf{x}}{{\rm d}p}$ is the tangent vector.

 To write out the formulas, we will use the bracket notation
 \begin{equation*}
 [\mathbf{a}, \mathbf{b} ]=\det (\mathbf{a}, \mathbf{b} ), \qquad \mathbf{a}, \mathbf{b} \in \R^2,
 \end{equation*}
 to denote the cross product in the plane.
 Let $s$ be the centro-affine arc-length parameter, where
 $\mathbf{x}'=\fracz{{\rm d}\mathbf{x}}{{\rm d}s}$ is used to distinguish $\overset{\bigcdot}{\mathbf{x}}=\fracz{{\rm d}\mathbf{x}}{{\rm d}p}$.
We first note that a parametrized curve $\mathbf{x}(p)$ is regular if and only if it satisfies
 \begin{equation}\label{non-d}
 \bigl[ \mathbf{x} , \overset{\bigcdot}{\mathbf{x}} \bigr]
\neq 0, \qquad \bigl[ \overset{\bigcdot}{\mathbf{x}},\overset{\bigcdot\bigcdot}{\mathbf{x}}\bigr]\neq 0.
 \end{equation}
The required formulas are obtained by replacing the jet derivative coordinates $y_x,y_{xx}, \ldots $ by using the chain
rule to express $x$-derivatives in terms of $p$-derivatives, as given by
\begin{equation*}
 y_{nx} \longmapsto D_x^n y,
\end{equation*}
where
\begin{equation*}
 D_x = \frac 1{\overset{\bigcdot}{x}} \frac{{\rm d}}{{\rm d}p}
\end{equation*}
is the differentiation operator dual to the horizontal one-form ${\rm d}x = \overset{\bigcdot}x\, {\rm d}p$.
Thus,
\begin{equation*}
 y_x \longmapsto \frac{\overset{\bigcdot}y}{\overset{\bigcdot}x}, \qquad y_{xx} \longmapsto \frac{\overset{\bigcdot}x\overset{\bigcdot\bigcdot}y - \overset{\bigcdot\bigcdot}x\overset{\bigcdot}y}{\overset{\bigcdot}x{}^3},
\end{equation*}
and so on. Substituting into \eqref{aff_arcl} produces the general formula for the centro-affine arc-length element of a parametrized curve
 \begin{equation}\label{arc-p}
 {\rm d}s
 =
 \sqrt{
 \epsilon\frac{
 \bigl[\overset{\bigcdot}{\mathbf{x}},\overset{\bigcdot\bigcdot}{\mathbf{x}}\bigr] } {
 \bigl[\mathbf{x},\overset{\bigcdot}{\mathbf{x}}\bigr] }} {\rm d}p,
 \qquad \textrm{where} \quad
\epsilon=\operatorname{sign}
 \left(
 \frac{
 \bigl[\overset{\bigcdot}{\mathbf{x}},\overset{\bigcdot\bigcdot}{\mathbf{x}}\bigr] }
 { \bigl[\mathbf{x},\overset{\bigcdot}{\mathbf{x}}\bigr] }
 \right). \end{equation}
One easily verifies that $\epsilon$ is invariant under centro-affine transformations and repara\-metrizations, including those that are orientation reversing.
Similarly, its centro-affine curvature \eqref{aff_curv} is given by the general formula
 \begin{gather}
 \sqrt{
 \epsilon\frac{
 \bigl[\mathbf{x},\overset{\bigcdot}{\mathbf{x}}\bigr]
 }
 {
 \bigl[\overset{\bigcdot}{\mathbf{x}},\overset{\bigcdot\bigcdot}{\mathbf{x}}\bigr]
 }
 }
 \left(
 \frac{
 \bigl[\mathbf{x},\overset{\bigcdot\bigcdot}{\mathbf{x}}\bigr]
 }{
 \bigl[\mathbf{x},\overset{\bigcdot}{\mathbf{x}}\bigr]
 }
 +\frac{
 \bigl[\mathbf{x},\overset{\bigcdot\bigcdot}{\mathbf{x}}\bigr]\bigl[\overset{\bigcdot}{\mathbf{x}}, \overset{\bigcdot\bigcdot}{\mathbf{x}}\bigr]
 -\bigl[\mathbf{x},\overset{\bigcdot}{\mathbf{x}}\bigr] \bigl[\overset{\bigcdot}{\mathbf{x}},\overset{\bigcdot\bigcdot\bigcdot}{\mathbf{x}}\bigr]
 }{
 2\bigl[\overset{\bigcdot}{\mathbf{x}},\overset{\bigcdot\bigcdot}{\mathbf{x}}\bigr] \bigl[\mathbf{x},\overset{\bigcdot}{\mathbf{x}}\bigr] } \right).
 \label{k1_gp}
 \end{gather}
 It is easy to check that \eqref{k1_gp} is equivalent to \eqref{aff_curv}, up to a constant factor. For the sake of convenience,
 we choose \eqref{k1_gp} as the centro-affine curvature $\kappa$.
 In particular, parametrizing the curve by centro-affine arc-length $s$,
 that is,
\[\frac{\bigl[\mathbf{x}', \mathbf{x}''\bigr]}{[\mathbf{x}, \mathbf{x}'\bigr]}=\epsilon, \qquad \hbox{which implies} \quad
\frac{\bigl[\mathbf{x}', \mathbf{x}'''\bigr]\bigl[\mathbf{x}, \mathbf{x}'\bigr]-\bigl[\mathbf{x}', \mathbf{x}''\bigr]\bigl[\mathbf{x}, \mathbf{x}''\bigr]}{\bigl[\mathbf{x}, \mathbf{x}'\bigr]^2}=0,\]
 one obtains the formula \cite{gw,olv-5,pk,wil}:
 \begin{equation}\label{kappa_def}
 \kappa =\frac{\bigl[\mathbf{x}'', \mathbf{x}\bigr]}{\bigl[\mathbf{x}', \mathbf{x}\bigr]},
 \end{equation}
 and
 \begin{equation}\label{Basic_str}
 \mathbf{x}''=\kappa \mathbf{x}'-\epsilon\mathbf{x}.
 \end{equation}

 Furthermore, from equation (\ref{Basic_str}), we have
 \begin{equation}\label{kappa_ds}
 \mathbf{x}'''=\big(\kappa _s+\kappa^2-\epsilon\big)\mathbf{x}'-\epsilon\kappa \mathbf{x},
 \end{equation}
 and
 \begin{equation*}
 \mathbf{x}''''=\big(\kappa''+3\kappa \kappa _s-2\epsilon\kappa +\kappa^3\big)\mathbf{x}'-\epsilon\big(2\kappa _s+\kappa^2-\epsilon\big)\mathbf{x}.
 \end{equation*}

\begin{figure}[t] \centering
 \begin{tabular}{cc}
 \includegraphics[width=.395\textwidth]{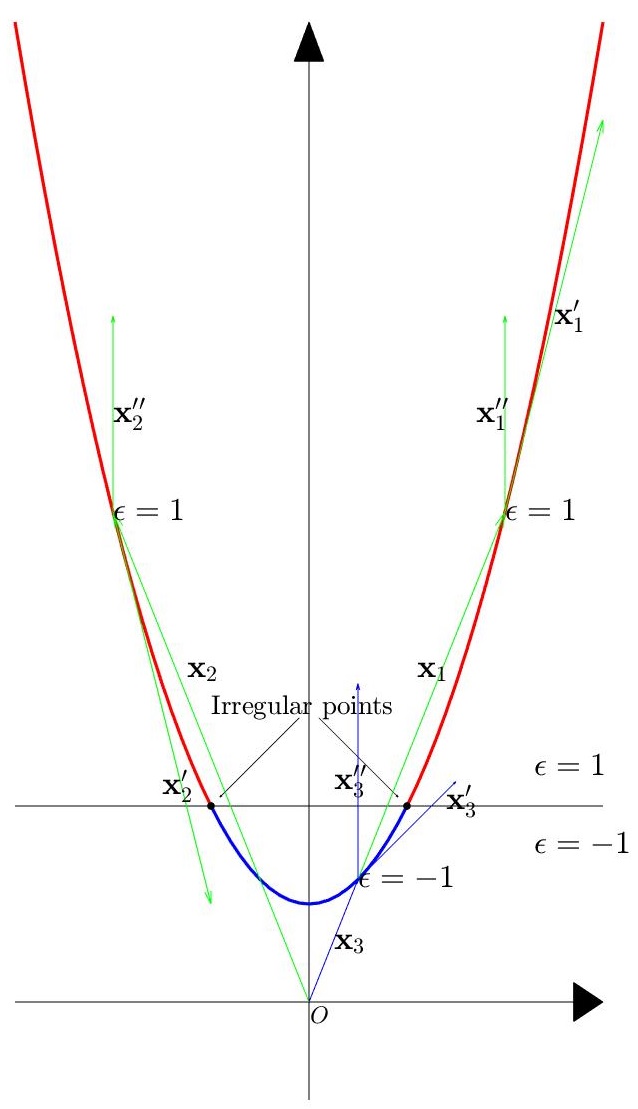}\qquad\qquad&\qquad\qquad\includegraphics[width=.295\textwidth]{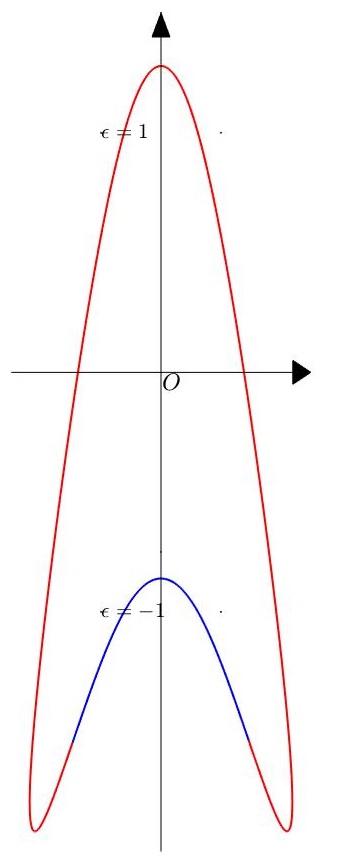}
 \end{tabular}
 \caption{\hbox{The geometrical properties of the centro-affine invariant $\epsilon$.}} \label{epsilon}
 \end{figure}

\begin{rem} For a regular curve $\mathbf{x}(p)$, $\epsilon\equiv1$ or $\epsilon\equiv-1$. The value of $\epsilon$ indicates that
 whether the vectors $\mathbf{x}$, $\mathbf{x}_{ss}$ lie on the same or opposite sides of the tangent vector $\mathbf{x}_s$.
 For example, in Fig.~\ref{epsilon}, $\epsilon=1$ on the red parts and $\epsilon=-1$ on the blue parts of the curve.
 The points in between the red and blue parts are irregular points.
 \end{rem}

 \begin{rem}
 If $\kappa (s)<0$, $s\in(a,b)$, let $\tilde{\kappa}(\tilde{s})=-\kappa (s)>0$, $\tilde{s}=-s\in(-b,-a)$.
 Then, by solving the differential equation
\[ \frac{{\rm d}^2\mathbf{x}(\tilde{s})}{{\rm d}\tilde{s}^2}-\tilde{\kappa}(\tilde{s})\frac{{\rm d}\mathbf{x}(\tilde{s})}{{\rm d}\tilde{s}}+\epsilon\mathbf{x}(\tilde{s})=0,\qquad \tilde{s}\in(-b,-a),
\]
 modeled on \eqref{Basic_str}, we obtain a curve $\tilde{\mathbf{C}}(\tilde{s})$ with the centro-affine curvature $\tilde{\kappa}(\tilde{s})>0$, $\tilde{s}\in(-b,-a)$.
 If we perform the orientation-reversing reparametrization
 $s=-\tilde{s}\in(a,b)$,
 then the resulting curve $\mathbf{C}(s)$ is equivalent to the curve $\tilde{\mathbf{C}}(\tilde{s}), \tilde{s}\in(-b,-a)$, and its centro-affine curvature satisfies $\kappa (s)<0$.
 \end{rem}

\section{Centro-affine planar curves with constant curvature}\label{sec-cai}

In any Klein geometry, submanifolds that have the property that all their differential invariants are constant play a particularly important role. The regular ones can all be algebraically characterized by the following theorem, originally due to Cartan \cite{Cartanrm,fo-2}. See also \cite{olv-sinmf,olv-sg} for further details.

\begin{thm} Let $G$ be a Lie group acting on an $m$-dimensional manifold $M$. Then, for $1 \leq p < m$, a regular $p$-dimensional submanifold $S \subset M$ has all constant differential invariants if and only if it is a subset of an orbit of a $($suitable$)$ $p$-dimensional subgroup $H \subset G$. In this case, $H$ can be identified with the local symmetry group of $S$.
 \end{thm}

The totally singular $p$-dimensional submanifolds, as defined above, are characterized by their admitting a (local) symmetry group of dimension $> p$. They can be characterized differentially by a Lie determinant condition~\cite{olv-sinmf}. The caveat ``suitable'' in the above result means that the subgroup $H$ has $p$-dimensional orbits and that such orbits are not totally singular. See~\cite{olv-sinmf} for an algebraic characterization.

Specializing to plane curves, we assume the transformation group is \emph{ordinary}, meaning that it acts transitively and does not pseudo-stabilize \cite{olv-eis}, which is the case for all fundamental geometric transformation groups including the centro-affine and centro-equi-affine groups. In this case, the differential invariant algebra is generated by invariantly differentiating a single differential invariant, which we identify as the $G$-invariant curvature. Thus, by the above results, a regular curve $C \subset X$ has constant curvature if and only if it forms part of the orbit of a suitable one-parameter subgroup $H \subset G$, and admits $H$ as its (local) symmetry group. Two such curves are \emph{equivalent}, meaning the $\widetilde C = g \cdot C$ for some $g \in G$ is their corresponding symmetry groups are related by the adjoint map $\widetilde H = \operatorname{Ad} g \cdot H = g \cdot H \cdot g^{-1}$. Thus one can use the methods of classification of optimal subalgebras of the Lie algebra of $G$ to classify constant curvature curves up to equivalence. On the other hand, totally singular curves are characterized by their admitting local symmetry groups of dimension $\geq 2$. All other curves have at most a discrete group of local symmetries.

Now let us apply these considerations to the centro-affine group. The first remark is that a~curve is totally singular if and only if it is a straight line segment. This follows easily from the prolonged action, or by calculation of the Lie determinant. For the group $\GLR2$, there are three inequivalent families of one-parameter subgroups, respectively generated by one of the Lie algebra elements
\begin{equation}\label{gl2ig}
\left(\begin{matrix}1 & 0 \\0 & \alpha \end{matrix}\right), \qquad
 \left(\begin{matrix}1 & 1 \\0 & 1 \end{matrix}\right), \qquad
 \left(\begin{matrix}\cos \varphi & -\sin \varphi \\ \sin \varphi & \cos \varphi \end{matrix}\right) ,
 \end{equation}
 where $\alpha \ne 0$ and $0 < \varphi < \pi$. In particular, when $\alpha = 1$, the corresponding orbit is a straight line and hence totally singular.
Exponentiating the infinitesimal generators~\eqref{gl2ig} to determine the corresponding one-parameter subgroups and orbits, we deduce the following result.

 \begin{thm}\label{thm-sg}
 Let $\mathbf{C}$ be a constant centro-affine curvature curve with $\epsilon=1$.
\begin{itemize}\itemsep=0pt
\item[$(i)$] If $\kappa>2$, then $\mathbf{C}$ is centro-affine equivalent to a curve
 \begin{equation*}
 y=x^\alpha, \qquad x>0, \quad 0 < \alpha < 1,
 \qquad \textrm{where} \quad \kappa=\frac{1}{\sqrt{\alpha}}+\sqrt{\alpha}.
 \end{equation*}
\item[$(ii)$] If $\kappa=2$, then $\mathbf{C}$ is centro-affine equivalent to the curve $y=x\log x$.
\item[$(iii)$]
If $\kappa = 0$, then $\mathbf{C}$ is centro-affine equivalent to the unit circle.
\item[$(iv)$]
If $0 < \kappa<2$, then $\mathbf{C}$ is centro-affine equivalent to the logarithmic spiral with polar coordinates
\begin{equation*}
 \rho=\exp(\theta \cot \varphi ),\qquad 0<\varphi <\pi/2, \qquad\textrm{where} \quad \kappa = 2 \cos \varphi .
\end{equation*}
\end{itemize}
 \end{thm}

 Note that the curves in Theorem~\ref{thm-sg} are generated from the subgroups in \eqref{gl2ig}, with $\varphi = \pi/2$ corresponding to the circle. The above values of the centro-affine curvature can be found either by direct computation or by applying the intrinsic algebraic method of~\cite{olv-max}.

 \begin{prop}A non-degenerate centro-affine curve has centro-affine curvature $\kappa=0$
 if and only if
 it is locally centro-affine equivalent to the unit hyperbola $y=x^{-1}$ or unit circle $x^2+y^2=1$.
 \end{prop}

\section{Centro-affine invariant geometric heat flow}\label{section4}

 Let us next investigate centro-affine invariant evolutionary processes.
 Consider a family of embedded smooth Jordan curves parametrized by $\mathbf{C}\colon S^1 \times I \rightarrow \R^2$, where $t \in I \subset \R$ can be viewed as the time parameter and $p \in S^1$ is a free parameter of each individual curve in the family.
 We assume that the curve family $\mathbf{C}(p,t)$ evolves according to the centro-affine invariant evolution equation
 \begin{equation}\label{ce-eq}
 \frac{\partial \mathbf{C}}{\partial t} = \beta(\kappa(s,t))\mathbf{C}_{ss}
 \end{equation}
 with the initial condition
 \begin{equation*}
 \mathbf{C}(s,0) = \mathbf{C}_0(s),
 \end{equation*}
 where $s$ is the corresponding centro-affine arc-length, $\kappa(s,t)$ is the centro-affine curvature, $\beta(\kappa)$~is a prescribed function of $\kappa$.

 Geometrically, equation (\ref{ce-eq}) means that any point of the curve moves with a velocity in the direction of the ``normal'' vector $\mathbf{C}_{ss}$, with speed proportional to a function of the centro-affine curvature of the curve at this point.
 These kinds of equations arise in differential geometry and a variety of applications due to their inherent invariance.
 On the other hand, because the arc length parameter $s$ changes with time, \eqref{ce-eq} is a non-linear evolutionary equation.

 In view of \eqref{arc-p}, let
 \begin{equation}\label{g-exp}
 g(p):=\sqrt{\epsilon \frac{[\mathbf{C}_p, \mathbf{C}_{pp}]}{[\mathbf{C}, \mathbf{C}_p]}}
 \end{equation}
 be the invariant centro-affine metric for the curve $C(p,t)$. The arc length parameter $s$ is obtained by integration:
 \begin{equation*}
 s(p)=\int^p_{p_0}g(\xi){\rm d}\xi.
 \end{equation*}
 In view of the elementary commutator relation
 $\displaystyle\frac{\partial}{\partial t}\frac{\partial }{\partial p} =\frac{\partial}{\partial p}\frac{\partial }{\partial t} $,
 we have
 \begin{equation*}
 \frac{\partial}{\partial t}\frac{\partial }{\partial s}
 =\frac{\partial}{\partial t}\left(\frac{1}{g}\frac{\partial}{\partial p}\right)
 =-\frac{g_t}{g}\frac{\partial}{\partial s}+\frac{\partial }{\partial s}\frac{\partial}{\partial t}.
 \end{equation*}

 Next, let us compute the {\it centro-affine metric evolution}. Firstly, using equation (\ref{g-exp}),
 \begin{equation*}
 \frac{\partial \big(g^2\big)}{\partial t}=\epsilon \frac{
 ( [\mathbf{C}_{pt}, \mathbf{C}_{pp} ]+ [\mathbf{C}_p, \mathbf{C}_{ppt} ] ) [\mathbf{C}, \mathbf{C}_p ]
 - [\mathbf{C}_{p}, \mathbf{C}_{pp} ] ( [\mathbf{C}_t, \mathbf{C}_p ]+ [\mathbf{C}, \mathbf{C}_{pt} ] )
 }{ [\mathbf{C}, \mathbf{C}_p ]^2 }.
 \end{equation*}
 Note that since the tangent $\mathbf{C}_s$ is not parallel to $\mathbf{C}$, then the right hand side of \eqref{ce-eq} can be expressed as the linear combination of $\mathbf{C}_s$ and $\mathbf{C}$, which means
 \begin{equation*}
 \beta\,\mathbf{C}_{ss}=W\mathbf{C}_s+U\mathbf{C},
 \qquad {\rm where} \quad W=\beta\kappa , \quad U=-\epsilon\beta.
 \end{equation*}
 By a direct computation, we obtain
 \begin{gather*}
 \mathbf{C}_p= g\mathbf{C}_s,\\
 \mathbf{C}_{pp}= gg_s\mathbf{C}_s+g^2\mathbf{C}_{ss}=\big(gg_s+g^2{\kappa }\big)\mathbf{C}_s-\epsilon g^2\mathbf{C},\\
 \mathbf{C}_{pt}= (W\mathbf{C}_s+U\mathbf{C})_p\\
\hphantom{\mathbf{C}_{pt}}{} = g(W_s\mathbf{C}_s+W\mathbf{C}_{ss}+U_s\mathbf{C}+U\mathbf{C}_s)=g(W_s+W{\kappa }+U)\mathbf{C}_s+g(U_s-\epsilon W)\mathbf{C},\\
 \mathbf{C}_{ppt}
= \big[gg_s(W_s+\kappa W+U)+g^2(W_{ss}+2\kappa W_s+(\kappa_s+\kappa^2-\epsilon) W+U_s+\kappa U)\big]\mathbf{C_s}\nonumber\\
\hphantom{\mathbf{C}_{ppt}=}{} +\big[gg_s(U_s-\epsilon W)-\epsilon g^2(2W_s+\kappa W-\epsilon U_{ss} +U)\big]\mathbf{C}.
 \end{gather*}
 Then
 \begin{gather*}
 [\mathbf{C}_{pt}, \mathbf{C}_{pp} ] =\epsilon g^2( g(W_s+U+\epsilon\kappa U_s)-g_s(W-\epsilon U_s)) [\mathbf{C}, \mathbf{C}_{s} ],\\
 [\mathbf{C}_p, \mathbf{C}_{ppt} ] =\epsilon g^2(g(2W_s+{\kappa } W-\epsilon U_{ss} +U)-\epsilon g_s(U_s-\epsilon W)) [\mathbf{C}, \mathbf{C}_{s} ],\\
 [\mathbf{C}_{p}, \mathbf{C}_{pp} ] =\epsilon g^3 [\mathbf{C}, \mathbf{C}_{s} ],\\
 [\mathbf{C}_t, \mathbf{C}_p ] =gU\left[\mathbf{C}, \mathbf{C}_s\right],\\
 [\mathbf{C}, \mathbf{C}_{pt} ] =g(W_s+W{\kappa }+U) [\mathbf{C}, \mathbf{C}_{s} ].
 \end{gather*}
 It follows that
 \begin{gather}\label{g2}
\frac 1g\frac{\partial g}{\partial t}=\;\frac 12 (2W_s+\epsilon {\kappa } U_s-\epsilon U_{ss})=\frac12 (\beta_{ss}+\kappa\beta_s)+ \kappa_s\beta.
 \end{gather}
 We now come to a crucial computation, namely the {\it centro-affine curvature evolution}
 \begin{gather*}
 \left.\frac{\partial {\kappa }}{\partial t}\right|_p
 =\frac{\partial}{\partial t}\frac{[\mathbf{C}, \mathbf{C}_{ss}]}{[\mathbf{C}, \mathbf{C}_s]} =\frac{
 ([\mathbf{C}_t, \mathbf{C}_{ss}]+[\mathbf{C}, \mathbf{C}_{sst}])[\mathbf{C}, \mathbf{C}_s]
 -[\mathbf{C}, \mathbf{C}_{ss}]([\mathbf{C}_t, \mathbf{C}_s]+[\mathbf{C}, \mathbf{C}_{st}])
 }{
[\mathbf{C}, \mathbf{C}_s]^2
 }.
 \end{gather*}
 A direct computation gives
 \begin{gather*}
 \mathbf{C}_{st} =\mathbf{C}_{ts}-\frac{g_t}{g}\mathbf{C}_s=(W\mathbf{C}_s+U\mathbf{C})_s-\frac{g_t}{g}\mathbf{C}_s=\left(W_s+{\kappa } W+U-\frac{g_t}{g}\right)\mathbf{C}_s+(U_s-\epsilon W)\mathbf{C},\\
 \mathbf{C}_{sst} =\mathbf{C}_{sts}-\frac{g_t}{g}\mathbf{C}_{ss}=\left[\left(W_s+{\kappa } W+U-\frac{g_t}{g}\right)\mathbf{C}_s+(U_s-\epsilon W)\mathbf{C}\right]_s-{\kappa }\frac{g_t}{g}\mathbf{C}_s+\epsilon\frac{g_t}{g}\mathbf{C}\nonumber\\
 \hphantom{\mathbf{C}_{sst}}{}
 =\left[W_{ss}+\big(\kappa_s+\kappa^2-\epsilon\big)W+2U_s+\kappa \left(2W_s+ U-2\frac{g_t}{g}\right)-\left(\frac{g_t}{g}\right)_s\right]\mathbf{C}_s\nonumber\\
\hphantom{\mathbf{C}_{sst}=}{} -\epsilon\left[-2W_s+{\kappa } W-\epsilon U_{ss}+U-2\frac{g_t}{g}\right]\mathbf{C}.
 \end{gather*}
 Thus we arrive at
 \begin{gather}\label{ka_ev}
 \frac{\partial {\kappa }}{\partial t}= {\kappa }_s W+2U_s-\frac{\epsilon}{2}\big({\kappa }_sU_s+{\kappa }^2U_s-U_{sss}\big)
 =\beta{\kappa }{\kappa }_s-2\epsilon\beta_s+\frac{1}{2}\big(\big(\kappa_s+{\kappa }^2\big)\beta_s-\beta_{sss}\big).\!\!\!
 \end{gather}

 We now focus on the case of $\beta \equiv 1$, namely the \emph{heat flow} in centro-affine geometry
 \begin{equation}\label{caheat}
 \frac{\partial \mathbf{C}}{\partial t} = \mathbf{C}_{ss}.
 \end{equation}
Equation (\ref{g2}) implies that
 \begin{equation}\label{gt-g}
 g_t=g\kappa _s.
 \end{equation}
 Consequently, in view of equation (\ref{ka_ev}), we see that the centro-affine curvature satisfies the first order inviscid Burgers' equation
 \begin{equation}\label{k1_t}
 \frac{\partial \kappa }{\partial t}=\kappa \kappa _s.
 \end{equation}

 To sum up, we arrive at the following results.
 \begin{thm}
 The centro-affine curve evolution process
 \begin{equation*}
 \frac{\partial\mathbf{C}}{\partial t}=\mathbf{C}_{ss},\qquad
\mathbf{C}(s,0)=\mathbf{C}_0,
 \end{equation*}
 is equivalent to the initial problem of the inviscid Burgers' equation:
 \begin{equation*}
\frac{\partial\kappa }{\partial t}=\kappa \frac{\partial\kappa }{\partial s},\qquad
\kappa (s,0)=\kappa _0(s),
 \end{equation*}
 where $\kappa _0(s)$ is the signed centro-affine curvature of the initial curve~$\mathbf{C}$.
 \end{thm}

 Comparing equation (\ref{gt-g}) with equation (\ref{k1_t}), one has
 \begin{equation*}
 \frac{\partial\kappa }{\partial t}-\frac{\kappa }{g}\frac{\partial g}{\partial t}=0.
 \end{equation*}
 Hence
 \begin{equation*}
 \frac{\partial (\fracz{\kappa }{g})}{\partial t}=\frac{1}{g}\left(\frac{\partial\kappa }{\partial t}-\frac{\kappa }{g}\frac{\partial g}{\partial t}\right)=0.
 \end{equation*}
 Thus we conclude:

\begin{cor}\label{e3.1}
 $\displaystyle \frac{\kappa }{g}(p,t)=\frac{\kappa }{g}(p,0)$ remains invariant for any $t>0$.
 \end{cor}

 Therefore, equation (\ref{k1_t}) can be written as
 \begin{equation}\label{kt-0}
 \frac{\partial \kappa }{\partial t}=\frac{\kappa }{g}(p,0)\frac{\partial\kappa }{\partial p}.
 \end{equation}
 In fact, if $\kappa _0(s) = \kappa (s,0)\equiv0$, by \eqref{kt-0}, the solution for \eqref{caheat} is $\mathbf{C}(\cdot,t)=\exp(-\epsilon t)\mathbf{C}_0$. We immediately deduce the following result.

 \begin{cor}
 If the initial centro-affine curvature $\kappa _0(s) = \kappa (s,0)\equiv0$,
 then at any time $t\geq0$, the curve $\mathbf{C}(t)$ is centro-affine equivalent to the initial curve $\mathbf{C}_0$.
 \end{cor}

 In the following, we assume $\kappa (p,0)\neq 0$. Solve the above partial differential equation \eqref{kt-0} by the method of characteristics yields
 \begin{equation*}
 \kappa (p,t)=\Phi(t+h(p)),
 \end{equation*}
 where $\Phi$ is any differentiable function of one variable and $\displaystyle h(p)=\int^{p}\frac{g}{\kappa }(\tilde{p},0)\,{\rm d}\tilde{p}$.
 Since $\kappa (p,0)\neq 0$, the ratio $\fracz{g(p,0)}{\kappa (p,0)}$ will be of one sign, and hence $h(p)$ defines a one-to-one map.
 Thus we conclude:
 \begin{cor}
 The curvature $\kappa (p,t)$ remains invariant on the curve $t+h(p)=C$, where~$C$ is constant.
 At the same time, $\kappa (p,t)=\kappa \big(h^{-1}(C),0\big)$,
 i.e., at any given time $t$, we have $\kappa (p,t)=\kappa \big(h^{-1}(t+h(p)),0\big)$.
 \end{cor}
 The flow (\ref{ce-eq}) can be written as
 \begin{equation}\label{Euclidean-flow}
 \frac{\partial \mathbf{C}}{\partial t}=\mathbf{C}_{ss}=\kappa \mathbf{C}_s-\epsilon\mathbf{C}.
 \end{equation}
 In view of Corollary~\ref{e3.1}, the heat flow~\eqref{Euclidean-flow} is equivalent to
 \begin{equation*}
 \epsilon \frac{\kappa}{g}(p,0)\mathbf{C}_p-\epsilon\mathbf{C}_t=\mathbf{C},
 \end{equation*}
 which can be solved easily, to get
 \begin{equation*}
 \mathbf{C}(p,t)=\exp(-\epsilon t)\mathbf{\Psi}\left(t+\int^p\frac{g}{\kappa}\,{\rm d}p\right),
 \end{equation*}
 where the vector function $\displaystyle \mathbf{\Psi}\left(\int^p\frac{g}{\kappa}\,{\rm d}p\right)=\mathbf{C}_0$, $\mathbf{C}_0=\mathbf{C}(p,0)$, $p\in(p_1,p_2)$, is the initial value of $\mathbf{C}(p,t)$, which implies
 \begin{equation}\label{Cau-C}
 \mathbf{\Psi}(\tilde{p})=\mathbf{C}_0\big(h^{-1}(\tilde{p})\big), \qquad \tilde{p}\in(h(p_1),h(p_2)).
 \end{equation}
 At any given time $\tau$, $\mathbf{\Psi}(\tilde{p}+\tau)=\mathbf{C}_0\big(h^{-1}(\tilde{p}+\tau)\big)$, where $\tilde{p}\in(h(p_1)-\tau, h(p_2)-\tau)$.

 On the other hand, by \eqref{g2}, the evolution of centro-affine arc-length $L$ for $\mathbf{C}(p,0)$ at $(p_1,p_2)$ is
 \begin{equation*}
 \frac{\partial L}{\partial t}=\frac{\partial}{\partial t}\int^{p_2}_{p_1} g \,{\rm d}p=\kappa(p_2)-\kappa(p_1).
 \end{equation*}
 Hence, we have
 \begin{cor}
 The solution to the initial problem of the heat flow \eqref{Euclidean-flow} with initial curve~$\mathbf{C}_0$, $p\in(p_1,p_2)$ and $\kappa\neq 0$ is given by $\displaystyle \mathbf{C}(p,t)=\exp(-\epsilon t)\mathbf{\Psi}\left(t+\int^p\frac{g}{\kappa}\,{\rm d}p\right)$ with $\mathbf{\Psi}$ satisfying~\eqref{Cau-C}.
 \end{cor}

 A basic fact on the theory of curve evolution~\cite{cz, eg} states that the geometric
 shape of the curve, sometimes referred to as the trace or the image of
 the planar curve, is affected only by the normal component of the flow field.
 The tangential component affects only the parametrization, and not the
 the curves' overall geometric shape. Thus, equation (\ref{Euclidean-flow}) is equivalent to
 \begin{equation}\label{flow-2}
 \frac{\partial \mathbf{C}}{\partial t}=-\epsilon\mathbf{C}.
 \end{equation}
 Solving the above equation, we have
 \begin{equation*}
 \mathbf{C}(p,t)=\exp(-\epsilon \,t)\mathbf{C}(p,0).
 \end{equation*}
 In this manner, we arrive at the long time behavior of the curves governed by the flow \eqref{Euclidean-flow}:

 \begin{prop} The curve family $\mathbf{C}(p,t)$ with $\epsilon=1$ will converge smoothly to the origin as $t \to \infty $.
 \end{prop}

\section{Edge matching of curve profiles in digital images}\label{section5}

 For many vision tasks, including 3D reconstruction, image alignment, and tracking, a key issue is finding correspondences between common objects in images.
 The SIFT and SURF algorithms are among the most widely applied to the identification of corresponding feature points.
In fact, when the images have less texture complexity and color diversity, it is not so easy to accurately extract and describe the feature points through application of the SIFT and SURF algorithms.
 Another drawback to these detectors is that they are only Euclidean-invariant.

 In general, the camera is often modeled as a projective transformation from scene coordinates to image coordinates.
 Especially, if a physical object has a smooth or piecewise smooth boundary, its images obtained by cameras in varying positions undergo apparent deformations, which are locally well approximated by affine transforms of the image plane.
 In consequence, the solid object recognition problem will lead back to the computation of affine invariant local image features.
 That is, during image matching, the Euclidean group (rotation, translation, reflection)
 should be extended to the equi-affine or full affine transformation group by including stretching and skewing transformations.

 In comparison with the above-mentioned methods, differential centro-affine-invariant detection has the following advantageous features:
 \begin{enumerate}\itemsep=0pt
 \item[(1)] It depends on the differential invariants of smooth curves for reducing errors.
 \item[(2)] Centro-affine differential invariants involve lower order derivatives of the curve paramet\-ri\-zation, and hence are more accurate and less error prone than fully affine differential invariants.
 \item[(3)] In some situations, centro-affine invariance is equivalent to fully affine invariance if we can find a pair of exact corresponding points (a point-correspondence) with respect to an affine transformation.
 \item[(4)] More precisely, for every closed curve, its barycenter can be temporarily served as that point-correspondence (or the local origin) only for its local centro-affine invariants.
 \item[(5)] The method relies solely on edge detection, which may well be adequate for untextured images.
 \end{enumerate}

\begin{figure}[hbtp] \centering
 \begin{tabular}{cccc}
 \includegraphics[width=.21\textwidth]{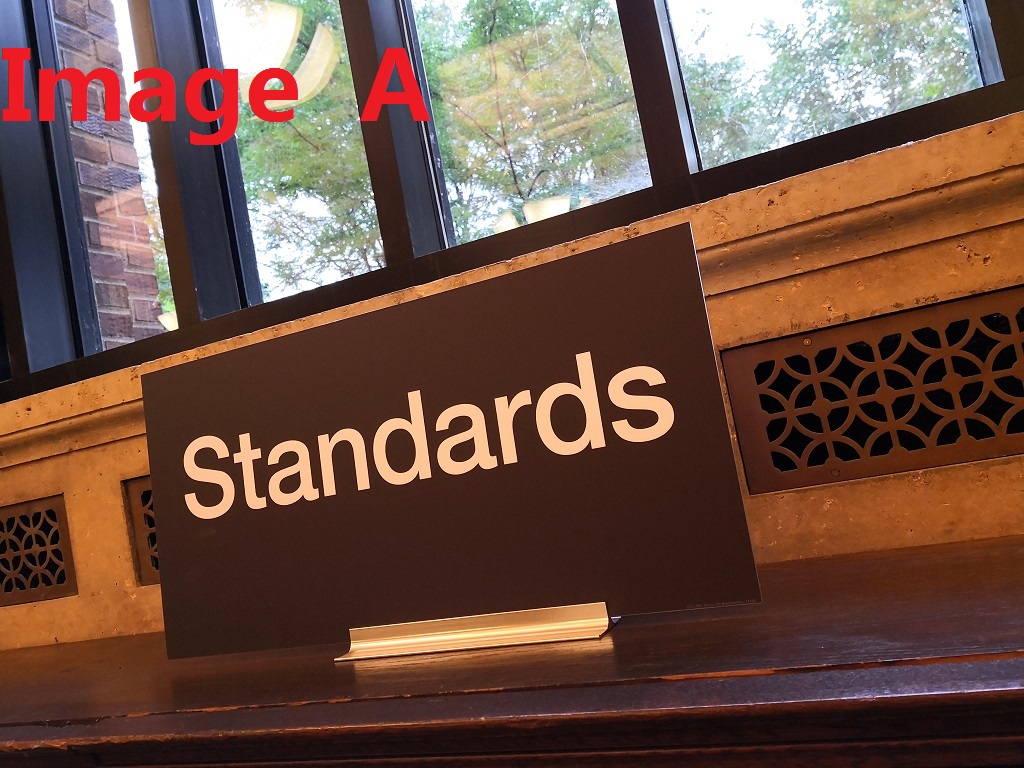}&\includegraphics[width=.22\textwidth]{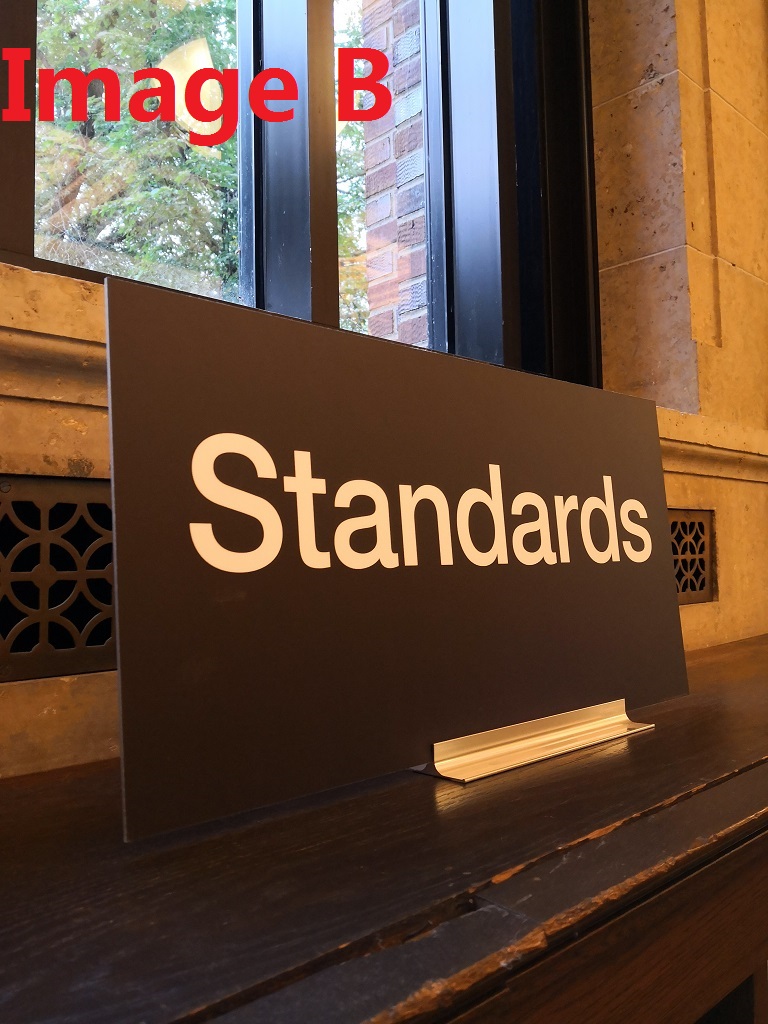}&\includegraphics[width=.2\textwidth]{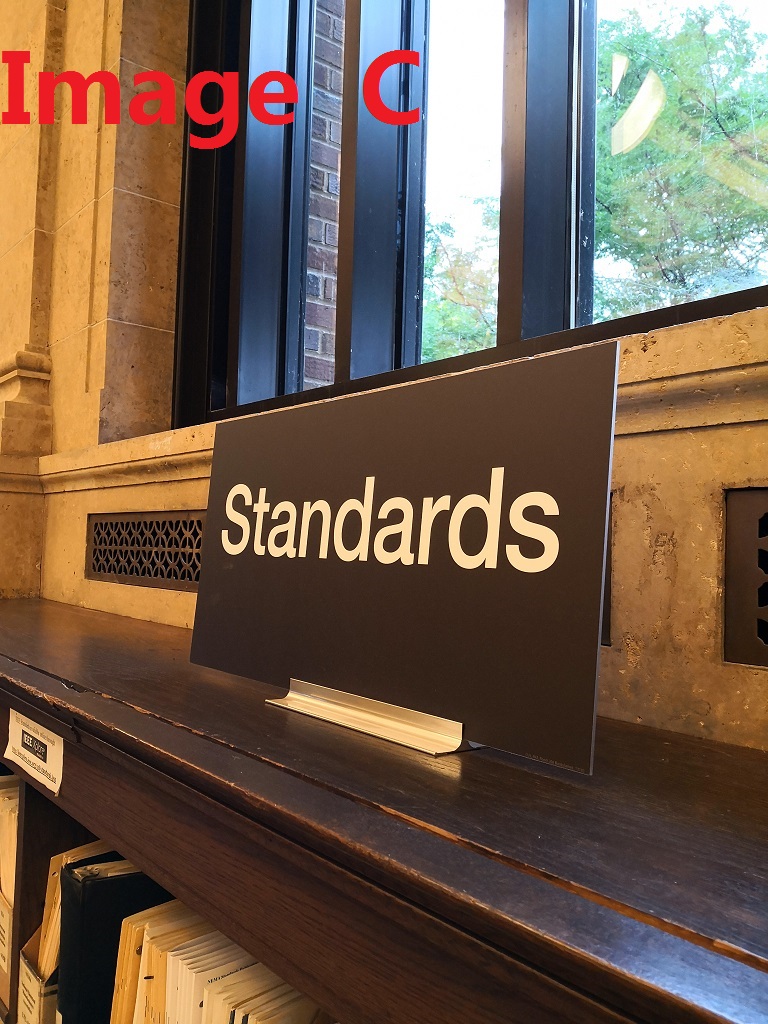}&\includegraphics[width=.2\textwidth]{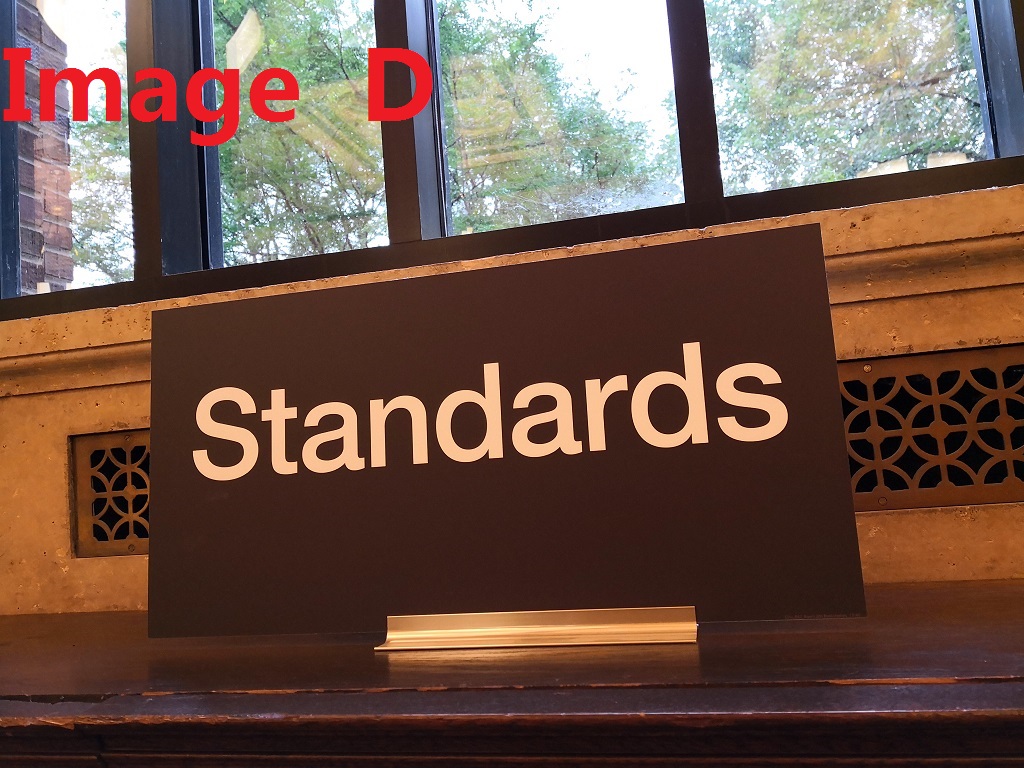}
 \end{tabular}
 \caption{The original images.} \label{or1-4}
 \end{figure}

In this section, we apply the centro-affine differential invariants defined in the previous sections to identify objects in different images.
 In Fig.~\ref{or1-4}, photos of a sign with the word ``Standards'' were taken by the camera at different angles and directions.
 We seek correspondences between the letters according to the centro-affine invariants given in equations (\ref{arc-p}) and (\ref{k1_gp}).

 \begin{figure}[hbtp] \centering
 \begin{tabular}{c}
 \includegraphics[width=.63\textwidth]{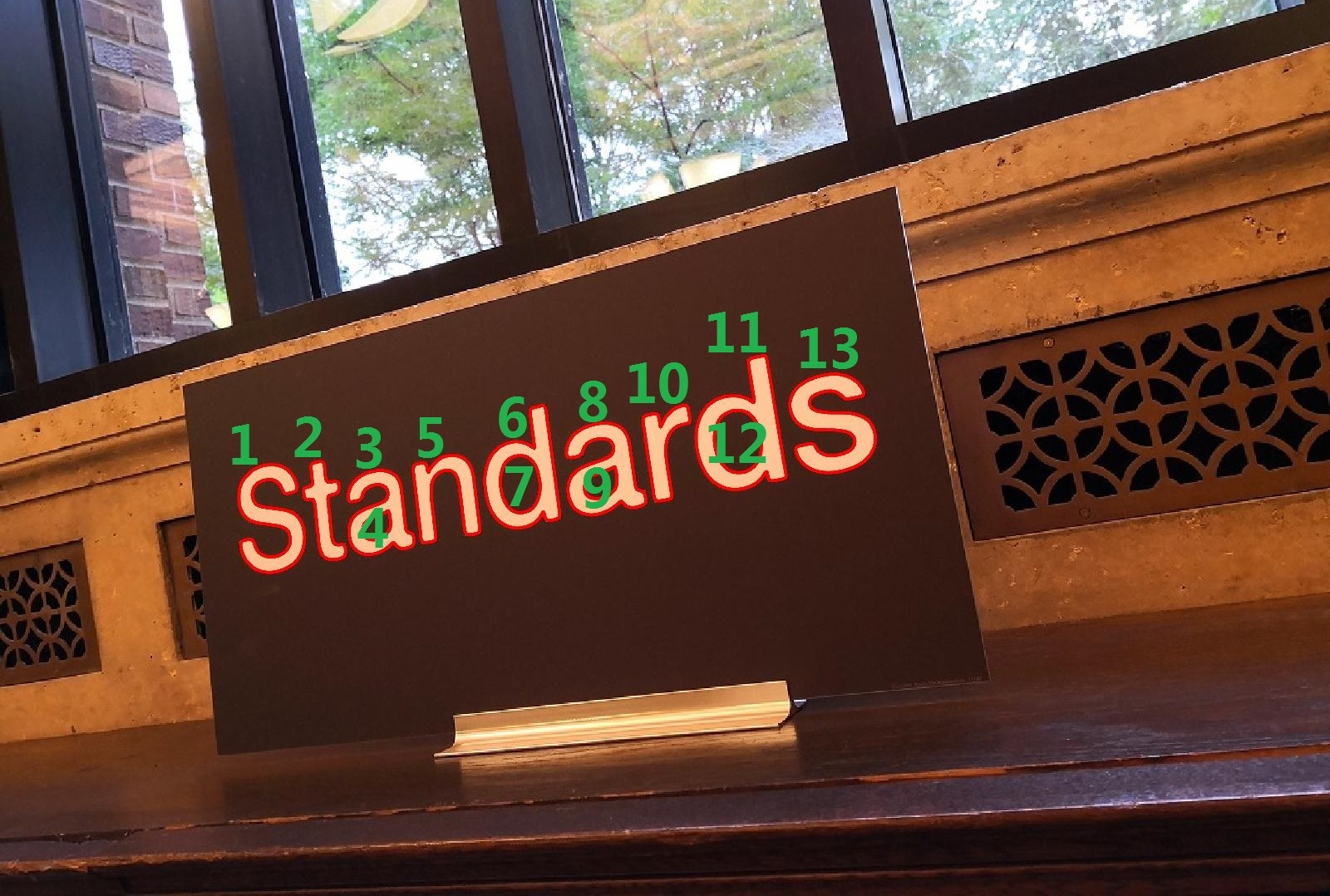}
 \end{tabular}
 \setlength{\belowdisplayskip}{3pt}
 \caption{The smooth edges generated by the B-spline curve fitting method and their ordinals, which will be used in Table~\ref{length}.} \label{orspline1-4}
 \end{figure}
 In order to describe more clearly and accurately this process, we now focus on the initial letter ``S'' as a main example; the others can be dealt with similarly.
 First of all, it is necessary to obtain the boundaries of the objects in the images.
 Accurate boundaries of
 the individual letters are segmented using level set method~\cite{lhd} or Canny edge detection~\cite{mos}. We shall apply the centro-affine arc-length,
 the integral invariant of centro-affine curvature, the area, the barycenter and corner points of every boundary as its descriptor,
 which can be employed to recognize the corresponding boundaries in different digital images.
 Corner point in a boundary represents critical information in describing
 object features, which is the local extreme point related to the Euclidean curvature.
 Then after the boundary extraction (segmentation), the letter ``S'' is located inside of its discrete boundary points.

 Notice that, for a given closed curve, its local centro-affine invariants may be obtained by placing the temporary origin at the barycenter of its closed boundary.
 To be precise, the smoothing algorithm~-- curve shortening flow~\cite{dde,gh} or Gaussian convolution for a curve~\cite{mos}~-- can play an auxiliary role in the pretreatment process for decreasing the jaggedness in the curve. For simplicity, we use Gaussian convolution for smoothing. The curve $\Gamma$ is first parameterized by the parameter $u$, so $\Gamma(u)=(x(u),y(u))$. An evolved version $\Gamma_\sigma$ of $\Gamma$ can then be computed, and is defined by $\Gamma_\sigma=(x(u,\sigma),y(u,\sigma))$,
 where $x(u,\sigma)=x(u) * g(u,\sigma)$, $y(u,\sigma)=y(u) * g(u,\sigma)$, with $*$ denoting the convolution operator, while $g(u, \sigma)$ denotes a~Gaussian of width~$\sigma$. For all boundaries, we apply k-means clustering~\cite{hw} to ensure that they are represented by the same number of points, which are used for further fitting purposes. By experimentation, we find that it is adequate to select~$85$ points for every boundary.
 To obtain the centro-affine invariants by~(\ref{arc-p}),~(\ref{k1_gp}), the discrete boundary points should generate a curve smooth enough for computing derivatives. B-spline curves are suitable for this role because of their affine invariance and smoothness \cite{chy, hc}. Thus, by a direct operation, we obtain the smooth boundaries as shown in Fig.~\ref{orspline1-4}.

\begin{figure}[hbtp]
 \centering
 \begin{tabular}{cc}
 \includegraphics[width=.41\textwidth]{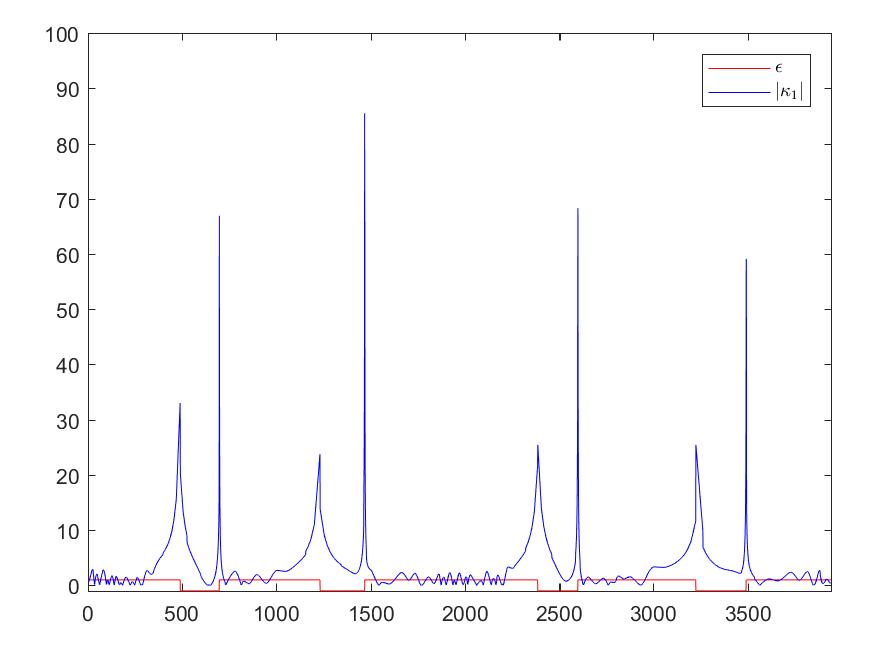}~&~\includegraphics[width=.41\textwidth]{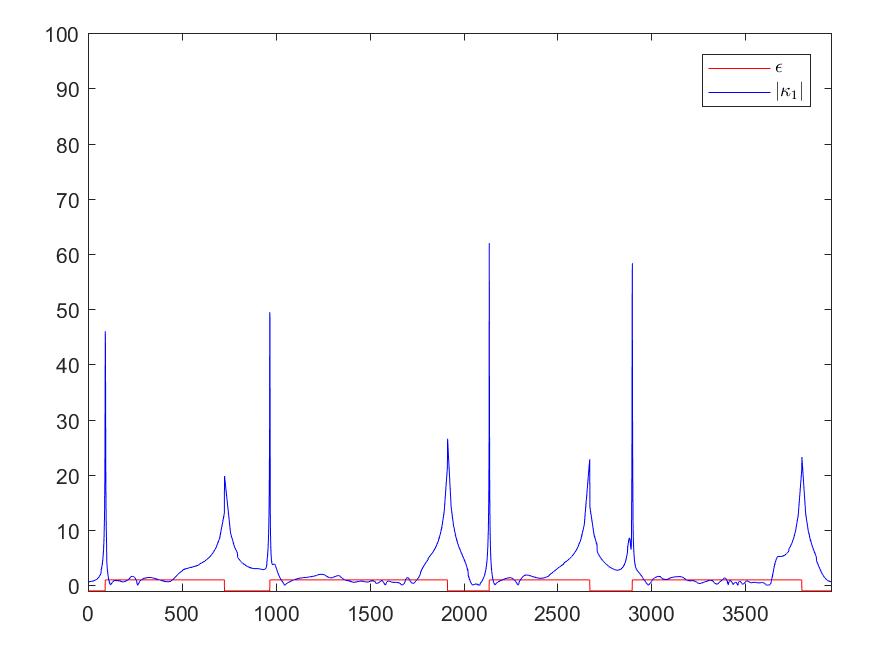}~\\
 \includegraphics[width=.41\textwidth]{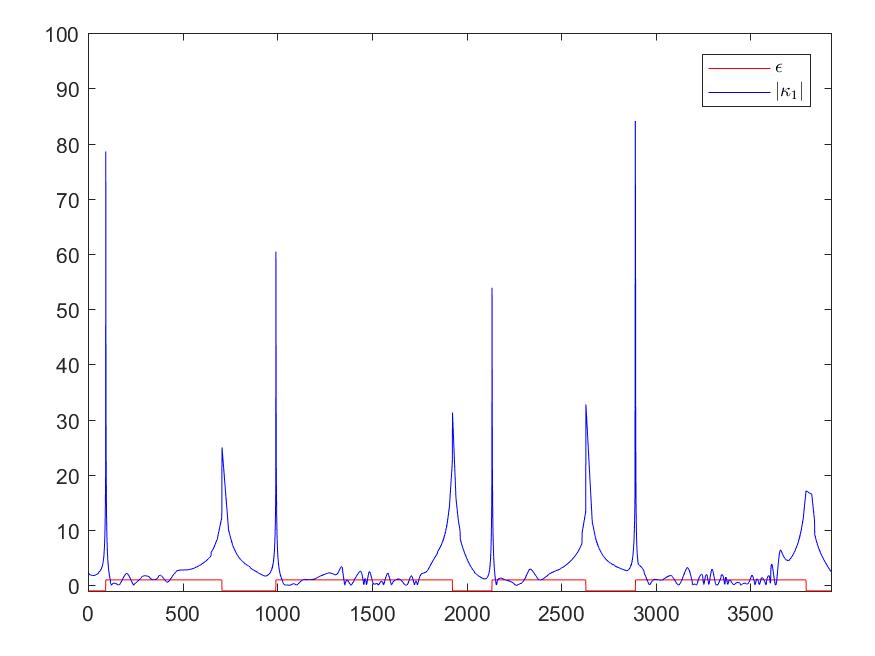}~&~\includegraphics[width=.41\textwidth]{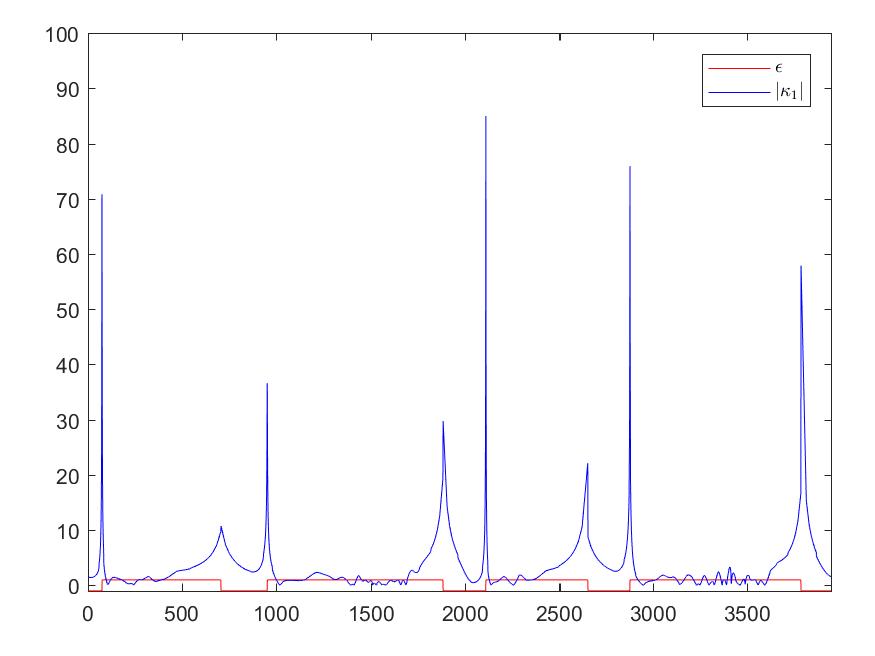}
 \end{tabular}
 \caption{The differential centro-affine invariants $\epsilon$ and $|\kappa |$ for the first letter ``S''. The horizontal axis represents the centro-affine arc-length up to a factor.} \label{curvature1-4}
 \end{figure}

If the boundary curve of the letter ``S'' is sufficiently smooth,
 it is easy to calculate its centro-affine invariants $\epsilon$ and centro-affine curvature $\kappa$ appearing in (\ref{arc-p}) and (\ref{k1_gp}), respectively.
 The final results are shown in Fig.~\ref{curvature1-4}.
 To reduce the disturbance caused by irregular points, we set threshold value of centro-affine curvatures $|\kappa |$ to $100$, deleting the points whose centro-affine curvatures $|\kappa |>100$.
 In view of the pictures in Fig.~\ref{curvature1-4}, we observe that they look similar modulo an overall translation in $s$.
 Now, it is crucial to find the corresponding points between them.

 We employ the ${\rm L}^2$ inner product
 \begin{eqnarray*}
 \ip{f(x)}{g(x)} =\int^b_af(x)g(x)\,{\rm d}x
 \end{eqnarray*}
 between functions $f(x)$ and $g(x)$ over the interval $[a,b]$. The associated norm{\samepage
\[ \operatorname{dist}(f(x),g(x)) = \nnorm{f(x)-g(x)}=\sqrt{\int^b_a(f(x)-g(x))^2\,{\rm d}x},\]
will be used as a measure of their distance.}

Now we can apply this concept to compare two discrete curves.
 Let us assume that, the curve~$A$ consists of $N$ points, denoted by
\[
 (x_{a_1}, y_{a_1})^{\rm T}, \ (x_{a_2}, y_{a_2})^{\rm T}, \ \ldots, \ (x_{a_N}, y_{a_N})^{\rm T}.
\]
 In this way, we denote the curve $B$ consists of $M$ points by
\[
 (x_{b_1}, y_{b_1})^{\rm T}, \ (x_{b_2}, y_{b_2})^{\rm T}, \ \ldots, \ (x_{b_M}, y_{b_M})^{\rm T}.
\]
 Locally,
 we can always assume that $x_{a_1}<x_{a_2}<\cdots<x_{a_N}$ and $x_{b_1}<x_{b_2}<\cdots<x_{b_M}$.
 For the sake of convenience, we use the following notations to indicate the maps of curve~$A$ and curve~$B$,
 that is,
\[
 y_{a_i}=f(x_{a_i}), \qquad y_{b_j}=g(x_{b_j}),\qquad i=1,2,\ldots,N,\quad j=1,2,\ldots,M.
\]
 Two new sets of $x$-coordinates are generated by trimming the original $\{x_{a_i}\}$ and $\{x_{b_i}\}$
 to only those values such that every $x'_{a_i}$ in the new set has two neighbours from $\{x_{b_i}\}$
 and every $x'_{b_i}$ in the new set has two neighbours from $\{x_{a_i}\}$.
 At the same time,
 we have $x'_{a_i}>x'_{b_{i-1}}$ and $x'_{b_i}>x'_{a_{i-1}}$.
 Then,
 we need to generate a common set of points for both curves with $x$-coordinates from the following set.
\[
 \{x_1,x_2,\ldots,x_K\}=\{x'_{a_1},x'_{a_2},\ldots,x'_{a_{N'}}\}\cup\{x'_{b_1},x'_{b_2},\ldots,x'_{b_{M'}}\},
\]
 where $x_1<x_2<\cdots<x_K$ and $\max\{M',N'\}\leq K\leq M'+N'$.
 The next step is to define the maps for the set $\{x_1,x_2,\ldots,x_K\}$,
 and here we take
 \begin{equation*}
 f(x_i)=
 \begin{cases}\displaystyle
 \frac{(x_i-x'_{a_{l-1}})f(x'_{a_l})+(x'_{a_l}-x_i)f(x'_{a_{l-1}})}{x'_{a_l}-x'_{a_{l-1}}}, & x'_{a_{l-1}}\leq x_i\leq x'_{a_l}, \quad 1<l\leq N', \\
 g(x_i),& x_i<x'_{a_1}, \\
 g(x_i),& x_i>x'_{a_{N'}},
 \end{cases}
 \end{equation*}
 and
 \begin{equation*}
 g(x_i)=
 \begin{cases}\displaystyle
 \frac{(x_i-x'_{b_{l-1}})g(x'_{b_l})+(x'_{b_l}-x_i)g(x'_{b_{l-1}})}{x'_{b_l}-x'_{b_{l-1}}}, & x'_{b_{l-1}}\leq x_i\leq x'_{b_l}, \quad 1<l\leq M', \\
 f(x_i),& x_i<x'_{b_1}, \\
 f(x_i),& x_i>x'_{b_{M'}}.
 \end{cases}
 \end{equation*}
 The missing $y$-coordinates (if any) for each curve are obtained via interpolating neighboring points.
 Thus, we can calculate the difference between curve~$A$ and curve~$B$ by the normalized~${\rm L}^1$ or~${\rm L}^2$ distance:
 \begin{equation*}
 \mathrm{error}_1=\frac{\mathrm{dist}_1}{\displaystyle\frac{1}{2N}\sum^N_{i=1}|y_{a_i}|+\frac{1}{2M}\sum^M_{i=1}|y_{b_i}|} \qquad \mathrm{or} \qquad
 \mathrm{error}_2=\frac{\mathrm{dist}_2}{\displaystyle\frac{1}{2N}\sqrt{\sum^N_{i=1}y^2_{a_i}}+\frac{1}{2M}\sqrt{\sum^M_{i=1}y^2_{b_i}}},
 \end{equation*}
where
 \begin{gather*}
 \mathrm{dist}_1=\max\{|f(x_i)-g(x_i)|, i=1,2,\ldots,K\}\qquad \mathrm{and}\qquad \mathrm{dist}_2=\frac{1}{K}\sqrt{\sum^K_{i=1}|f(x_i)-g(x_i)|^2}.
 \end{gather*}

 \begin{figure}[t] \centering
 \begin{tabular}{ccc}
 \includegraphics[width=.29\textwidth]{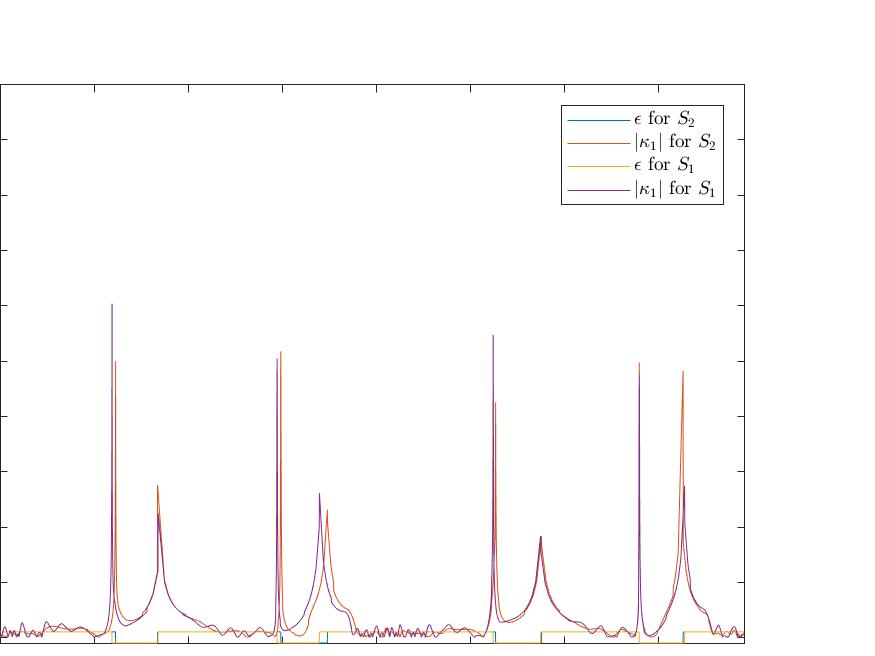}\ &\ \includegraphics[width=.29\textwidth]{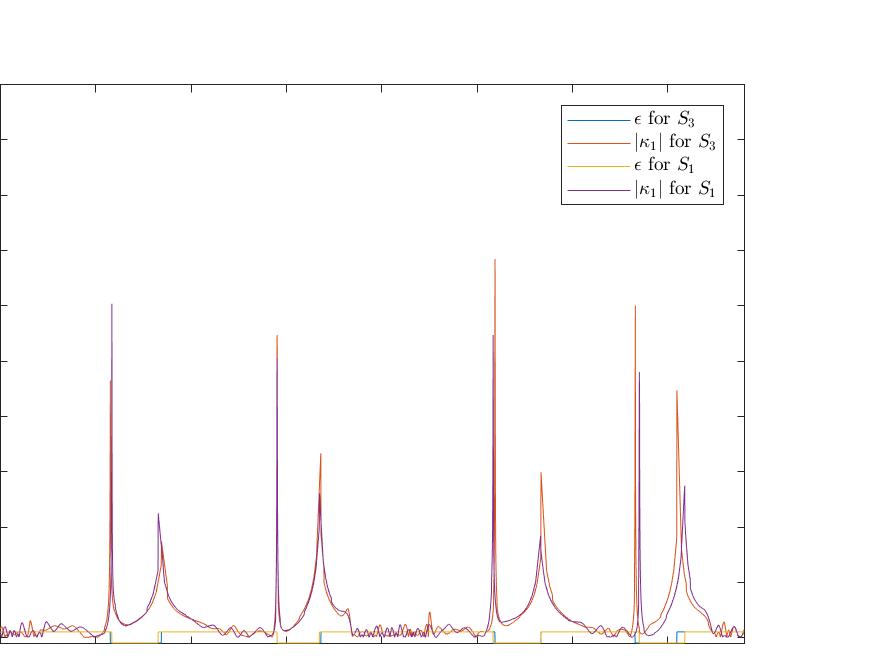}\ &\ \includegraphics[width=.29\textwidth]{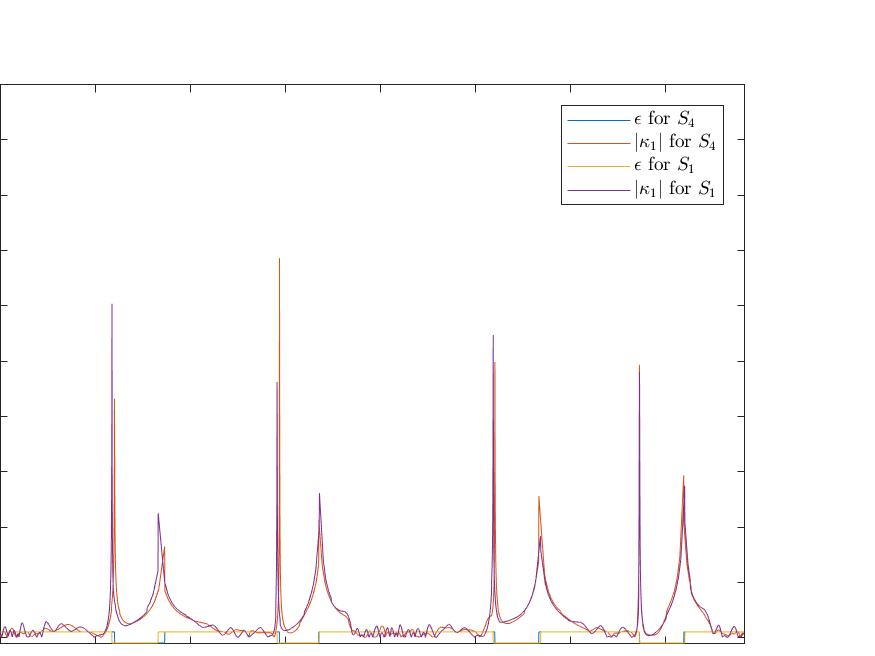}\\
 \includegraphics[width=.29\textwidth]{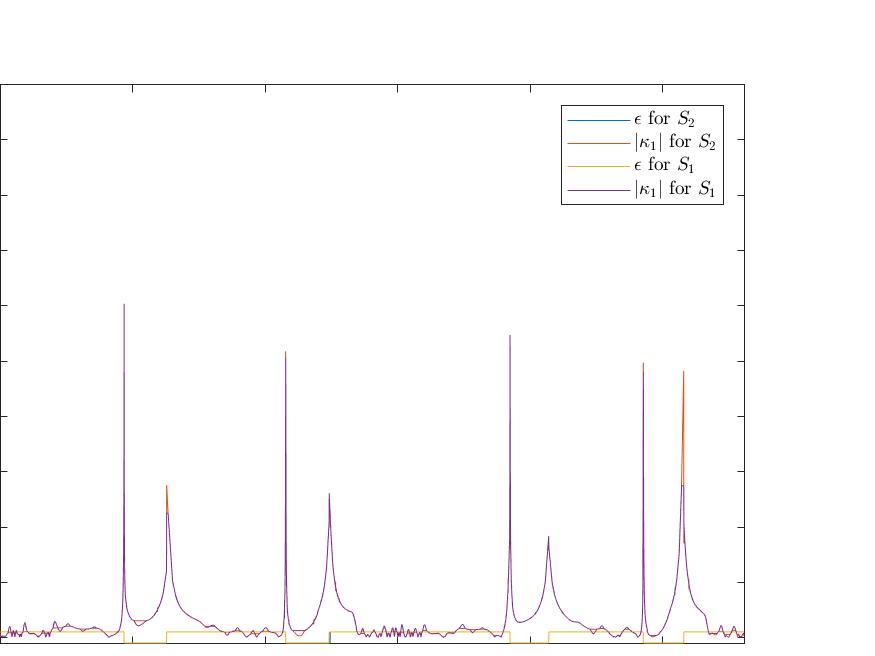}\ &\ \includegraphics[width=.29\textwidth]{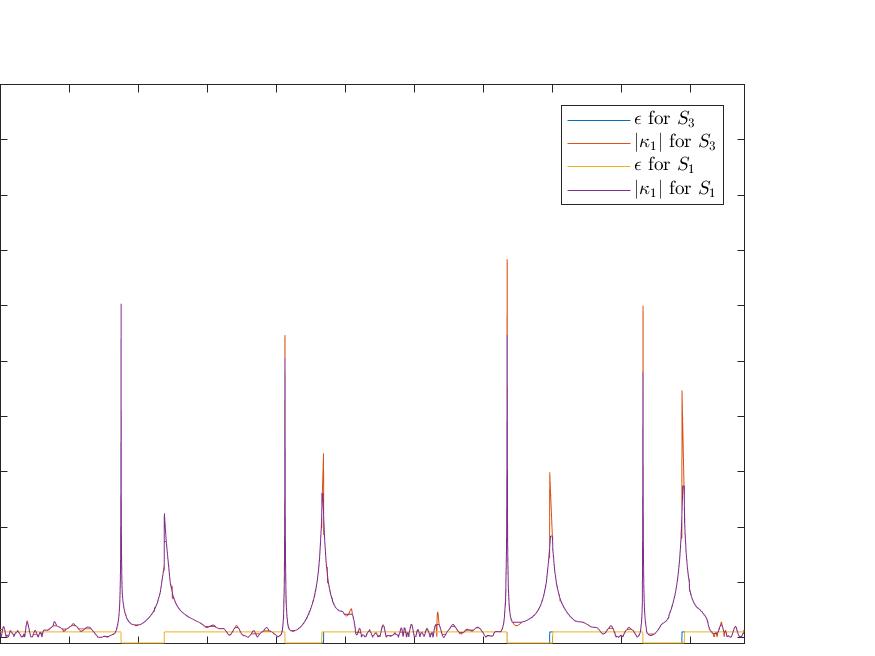}\ &\ \includegraphics[width=.29\textwidth]{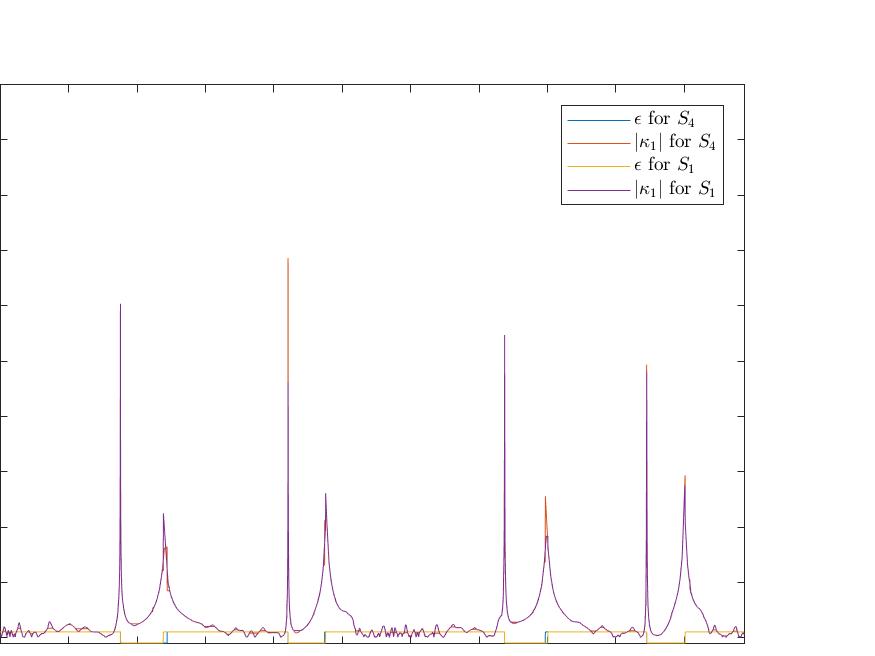}
 \end{tabular}

 \caption{Finding the corresponding relations by using the data in Fig.~\ref{curvature1-4}.} \label{compare}
 \end{figure}

\begin{table}[h!]\centering
 \footnotesize
 \caption{\label{length}The centro-affine arc-lengths for the smooth boundaries in Fig.~\ref{or1-4}, where the row labels and column labels correspond to the ordinals in Figs.~\ref{or1-4} and~\ref{orspline1-4}.} \vspace{1mm}

\begin{tabular}{|c|c|c|c|c|c|c|c|}
 \hline
 &1&2&3&4&5&6&7\\
 \hline
 Image A &21.5764&17.3975&16.0039&6.0827&13.1816&10.2993&6.1379\\
 \hline
 Image B &21.4343&18.1060&16.2536&6.0291&13.4718&9.9786&6.1225\\
 \hline
 Image C &21.3836&16.3080&15.0370&5.6557&12.2398&9.9733&6.1182\\
 \hline
 Image D &21.4965&17.0939&15.8373&5.9567&13.2712&11.4628&6.1752\\
 \hline
 ({\it continued})&8&9&10&11&12&13&\\
 \hline
 Image A &16.1042&6.0655&17.9754&10.5182&6.1844&20.9430&\\
 \hline
 Image B &16.3201&5.8780&12.8081&9.4719&6.0826&20.6943&\\
 \hline
 Image C &15.7962&5.8537&15.7287&11.1307&6.1087&20.9512&\\
 \hline
 Image D &16.7117&6.0081&18.1944&10.3007&6.1984&20.7571&\\
 \hline

 \end{tabular}
 \end{table}

 \begin{table}[h!]
 \centering\footnotesize
 \caption{\label{cor}The correlation coefficents for Table \ref{length}.}\vspace{1mm}

 \begin{tabular}{|c|c|c|c|c|}
 \hline
 &Image A&Image B&Image C&Image D\\
 \hline
 Image A &1&0.9672&0.9930&0.9976\\
 \hline
 Image B &0.9672&1&0.9767&0.9624\\
 \hline
 Image C &0.9930&0.9767&1&0.9897\\
 \hline
 Image D &0.9976&0.9624&0.9897&1\\
 \hline
 \end{tabular}\vspace{-2mm}
 \end{table}

 By using the above method, under translations or orientation reversals of data sets,
 we can find the minimum errors about $\epsilon$ and centro-affine curvatures between the first graph in the Fig.~\ref{curvature1-4} and the remaining three.
 In the first row of Fig.~\ref{compare}, the comparison results are shown together. We apply the dynamic time warping (DTW) algorithm~\cite{sc} to find the optimal alignment between the two sets of data points, which we view as time series, as shown in the second row of Fig.~\ref{compare}. In general, DTW is often used to determine similarity, classification, and corresponding regions between two time series.

\begin{figure}[t] \centering
 \begin{tabular}{ccc}
 \includegraphics[width=.22\textwidth]{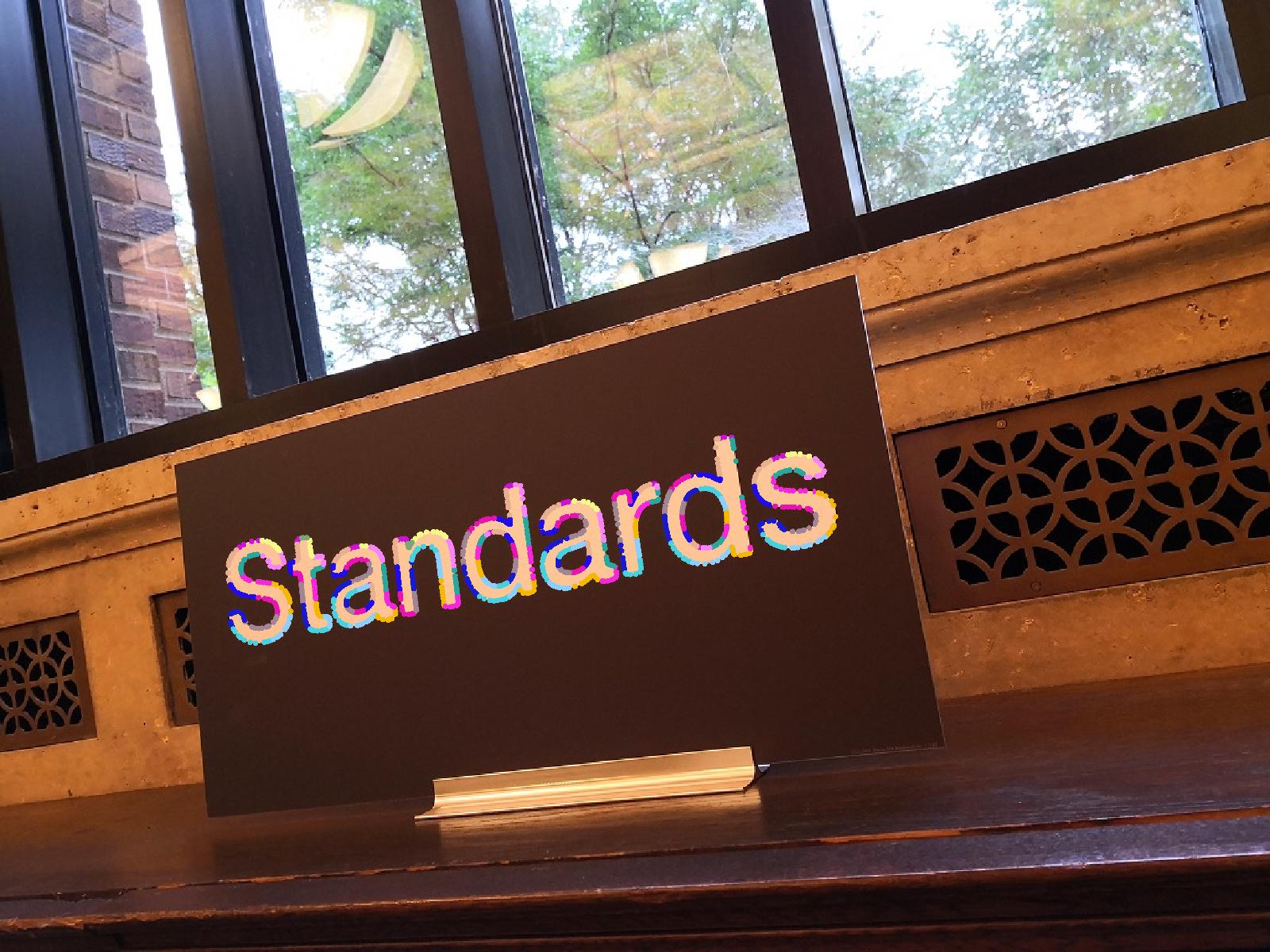}~~~\qquad&\qquad~~~\includegraphics[width=.22\textwidth]{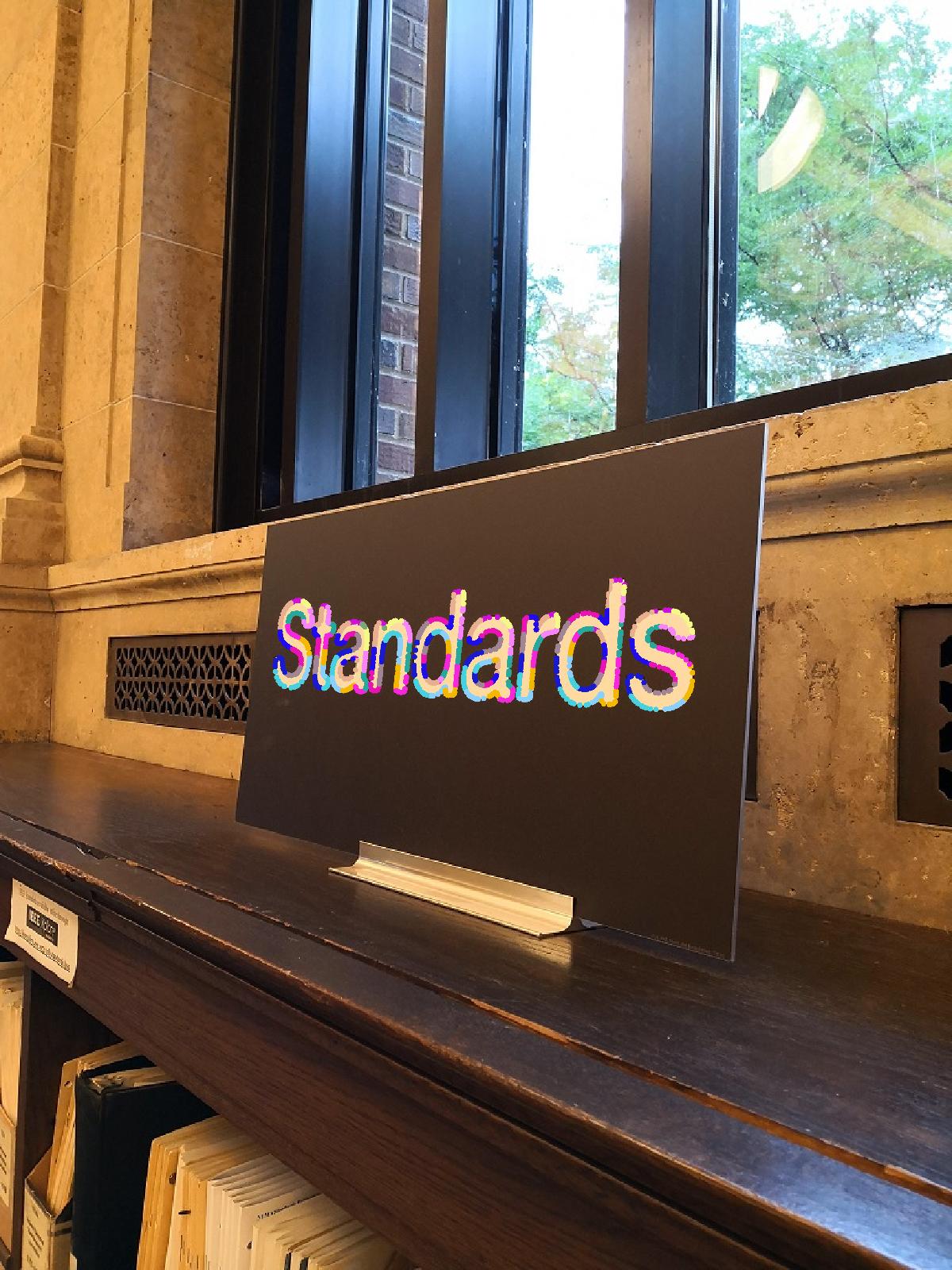}~~~\qquad&\qquad~~~\includegraphics[width=.22\textwidth]{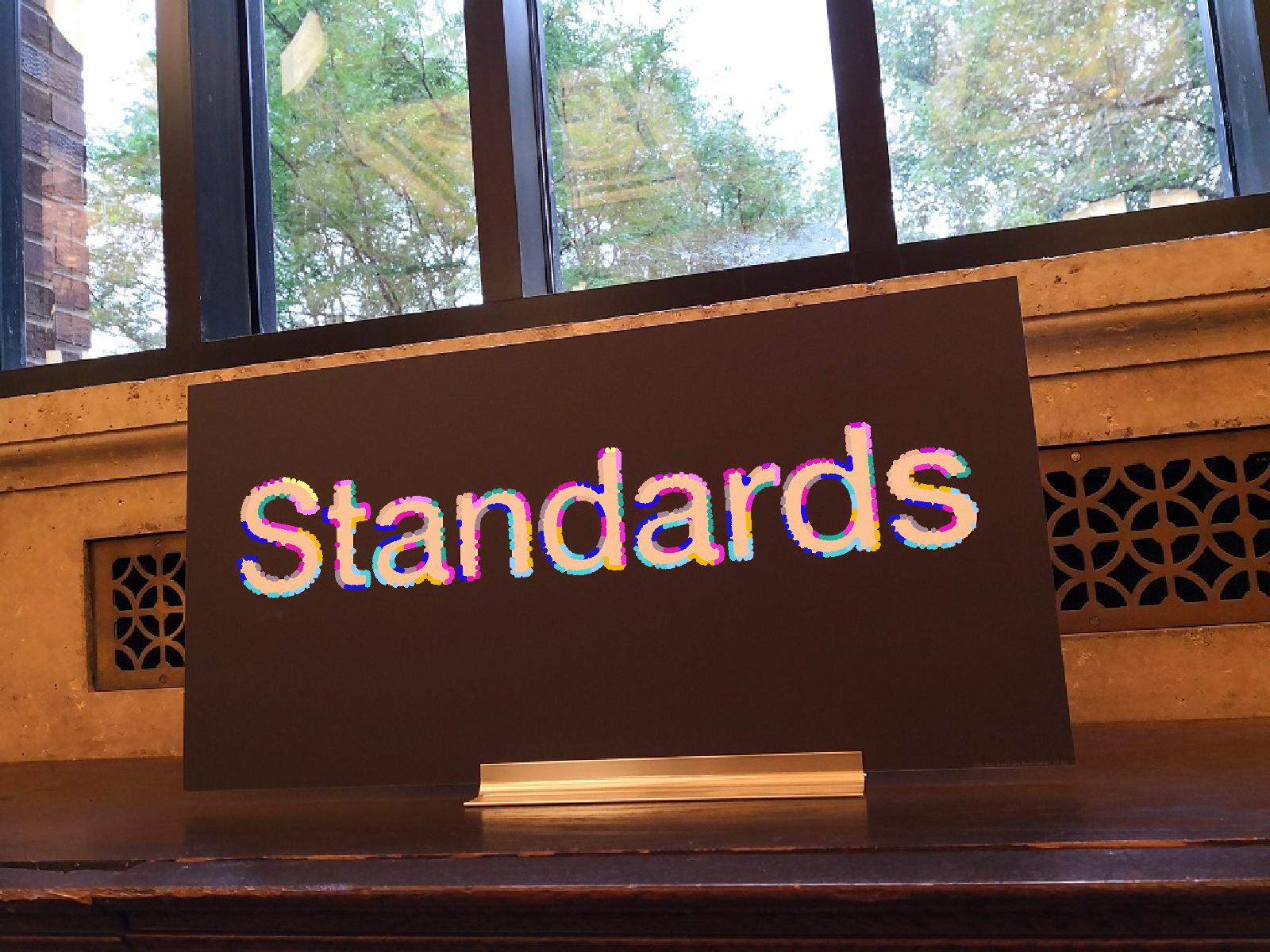}\\
 \includegraphics[width=.22\textwidth]{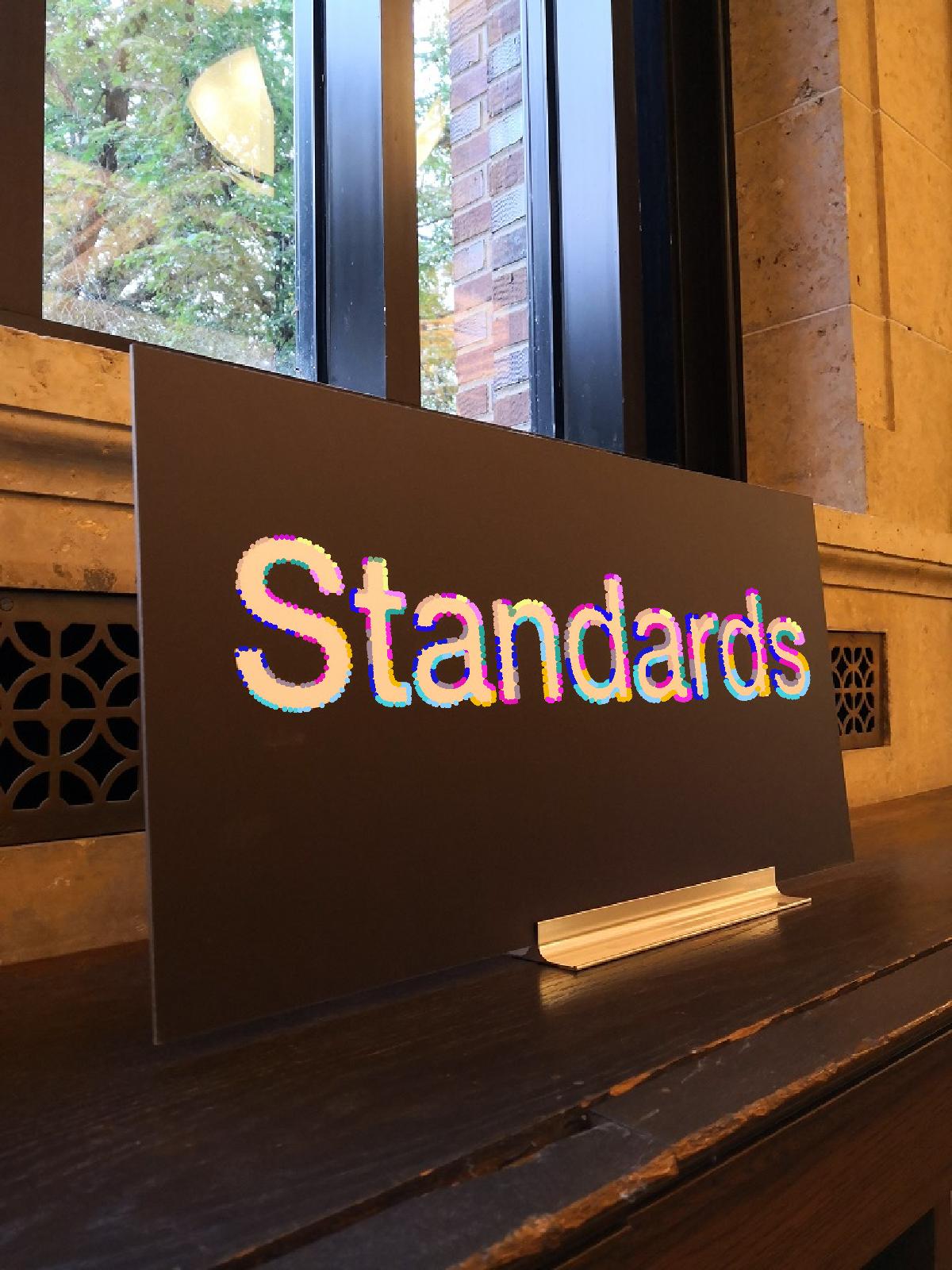}~~~\qquad&\qquad~~~\includegraphics[width=.22\textwidth]{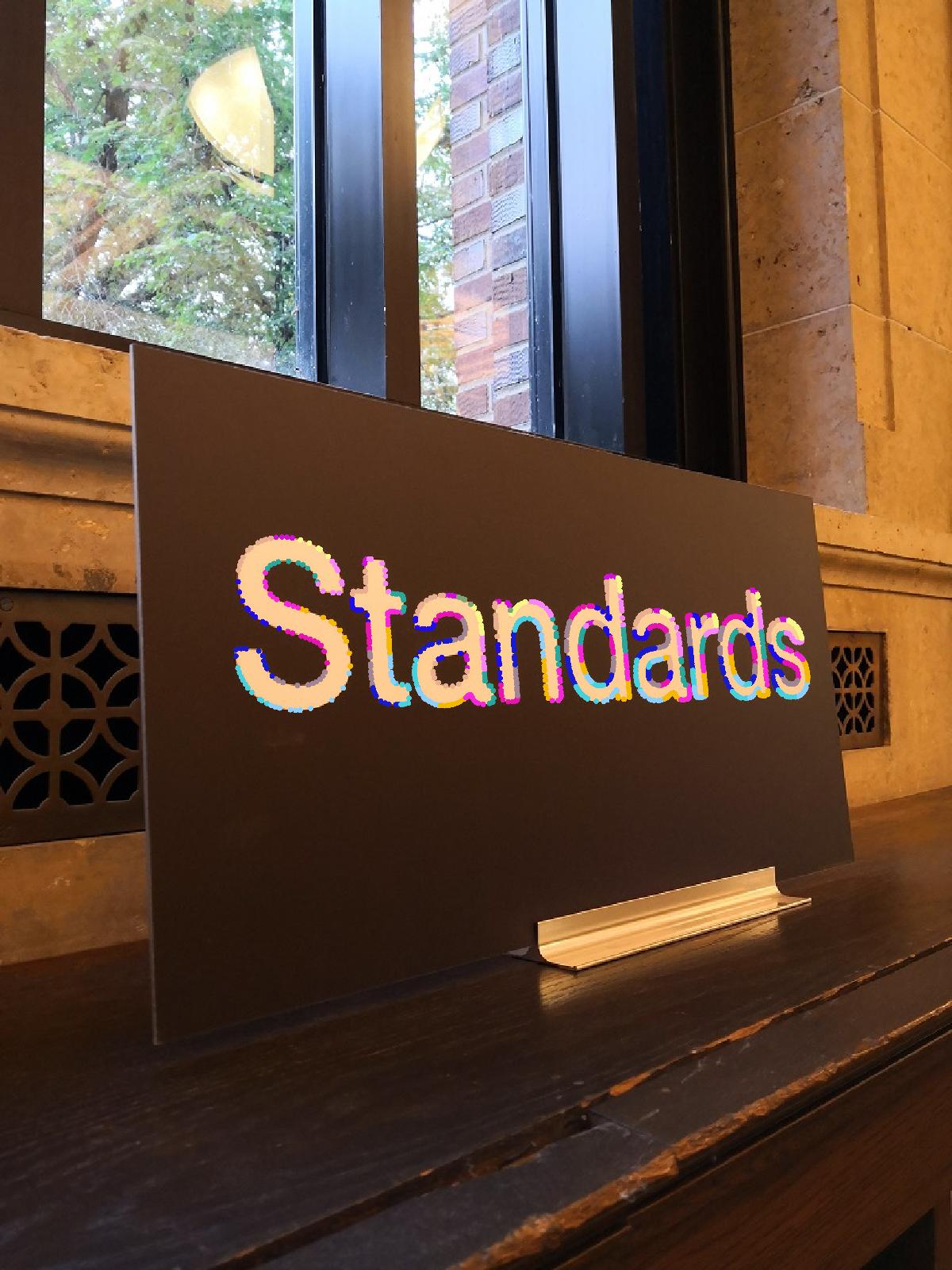}~~~\qquad&\qquad~~~\includegraphics[width=.22\textwidth]{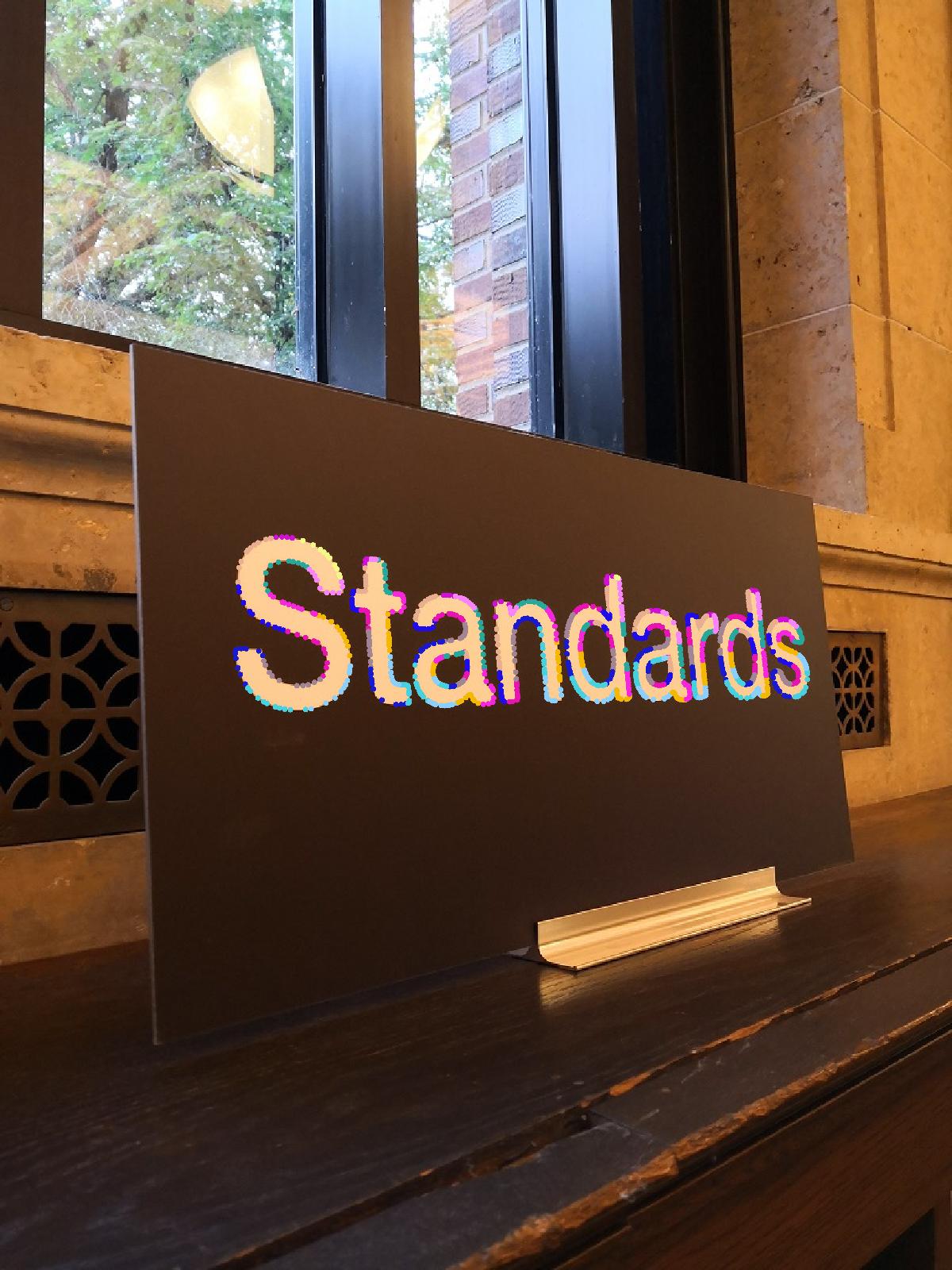}
 \end{tabular}
 \caption{The corresponding points between the two images of every column obtained by using the centro-affine invariants method.}
 \label{compareFig}
 \end{figure}

\begin{figure}[t] \centering
 \begin{tabular}{ccc}
 \includegraphics[width=0.24\textwidth,clip]{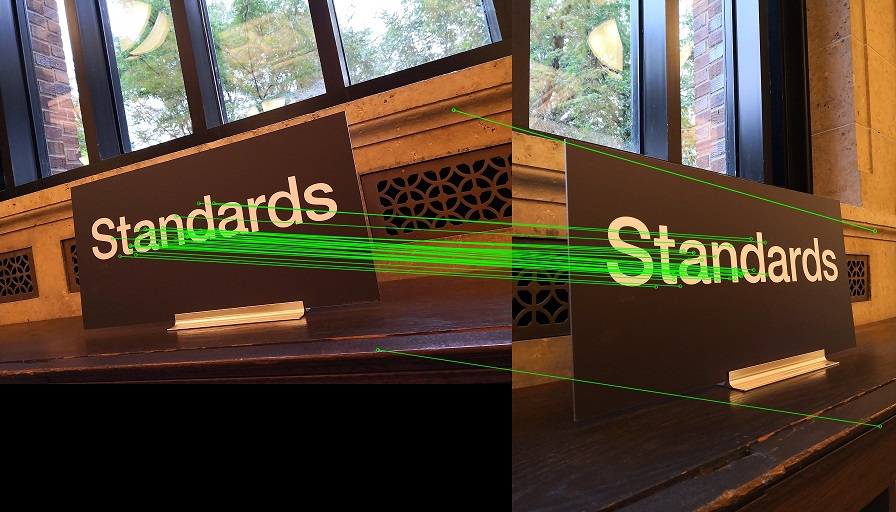}&\includegraphics[width=0.24\textwidth,clip]{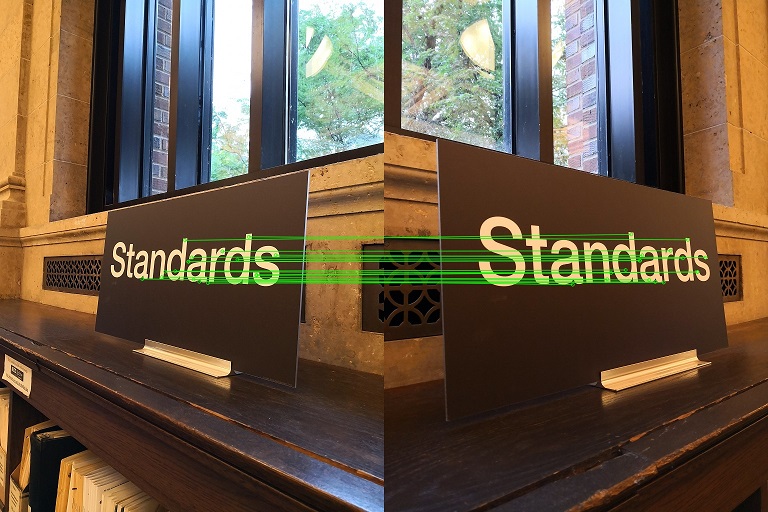}&\includegraphics[width=0.24\textwidth,clip]{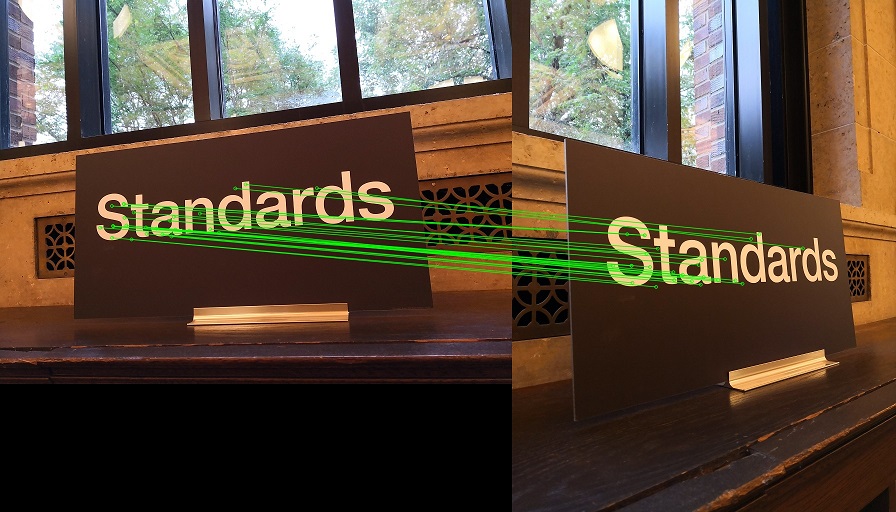}\\
 \includegraphics[width=0.24\textwidth,clip]{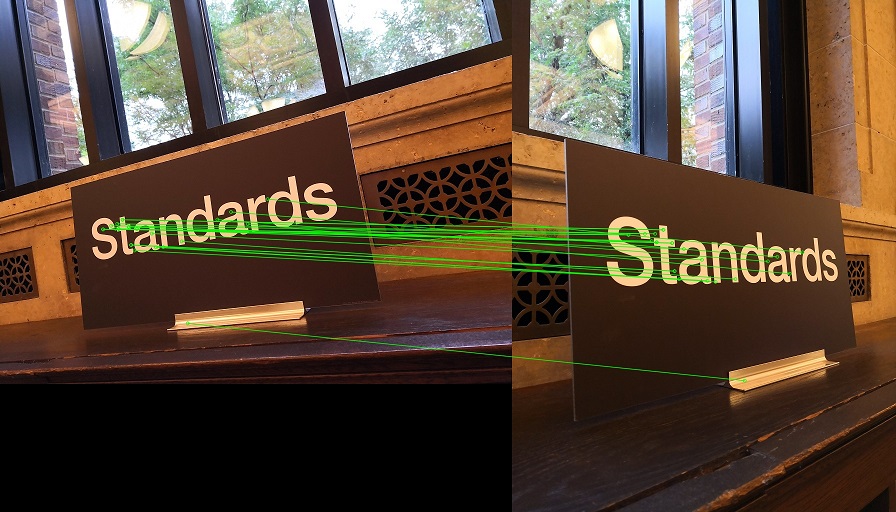}&\includegraphics[width=0.24\textwidth,clip]{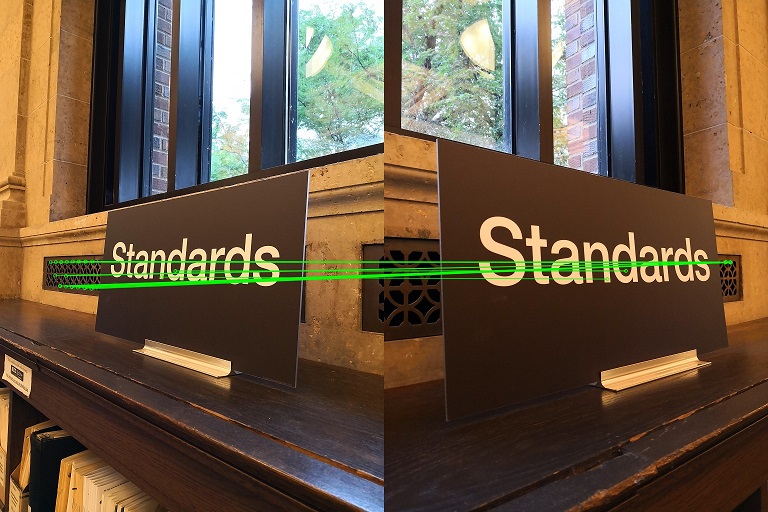}&\includegraphics[width=0.24\textwidth,clip]{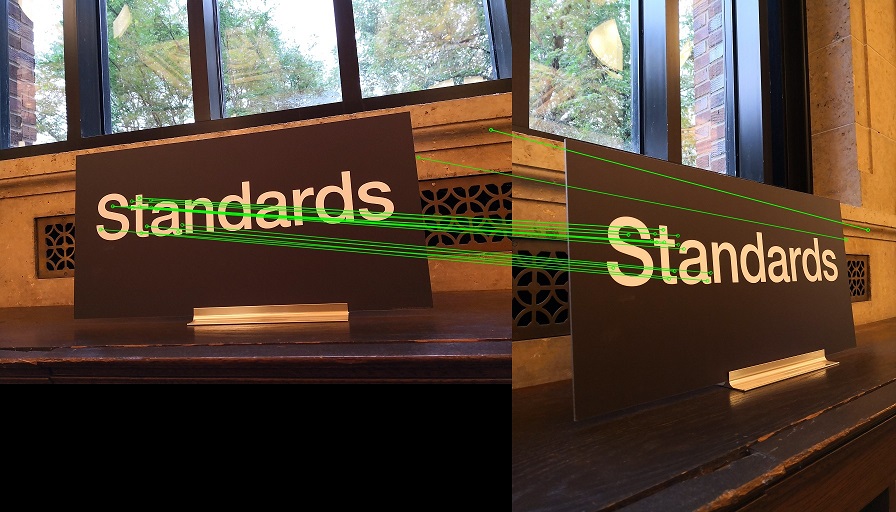}\\
 \includegraphics[width=0.24\textwidth,clip]{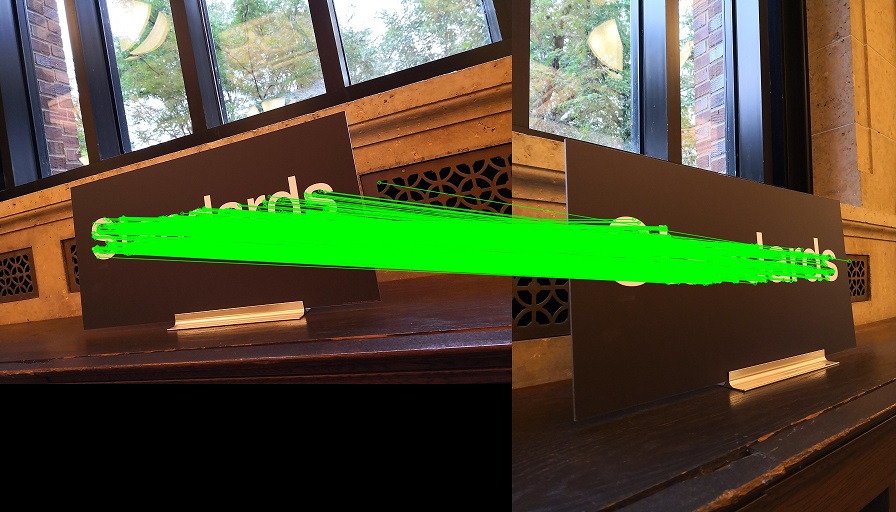}&\includegraphics[width=0.24\textwidth,clip]{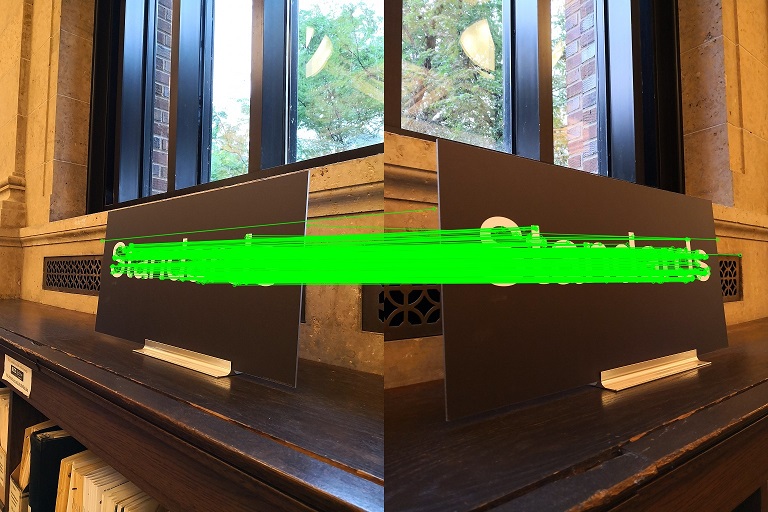}&\includegraphics[width=0.24\textwidth,clip]{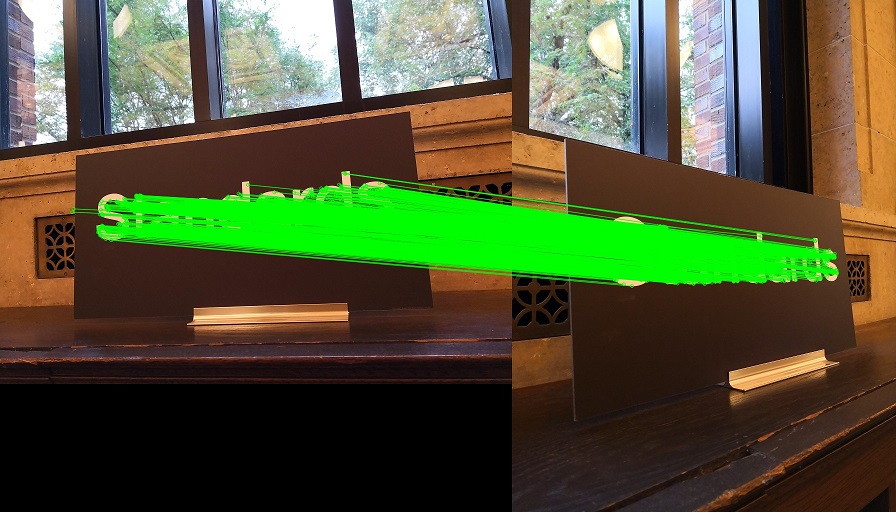}
 \end{tabular}
 \caption{The corresponding points obtained by the SURF method (the first row), the SIFT method (the second row) and the ASIFT method (the third row).} \label{siftcps}
 \end{figure}

 Finally, applying the same algorithm to the other letters in the label, the corresponding points can be found.
 The results are shown in Fig.~\ref{compareFig}, where, between two images of every column, the correspondences are indicated by the different colors.
 For comparative purposes, we also apply the SIFT, SURF and ASIFT methods to obtain their corresponding points,
 which are shown in Fig.~\ref{siftcps}, respectively.

 \begin{figure}[t] \centering
 \begin{tabular}{cccc}
 Centro-Affine&SURF&SIFT&ASIFT\\
 \includegraphics[width=0.18\textwidth,clip]{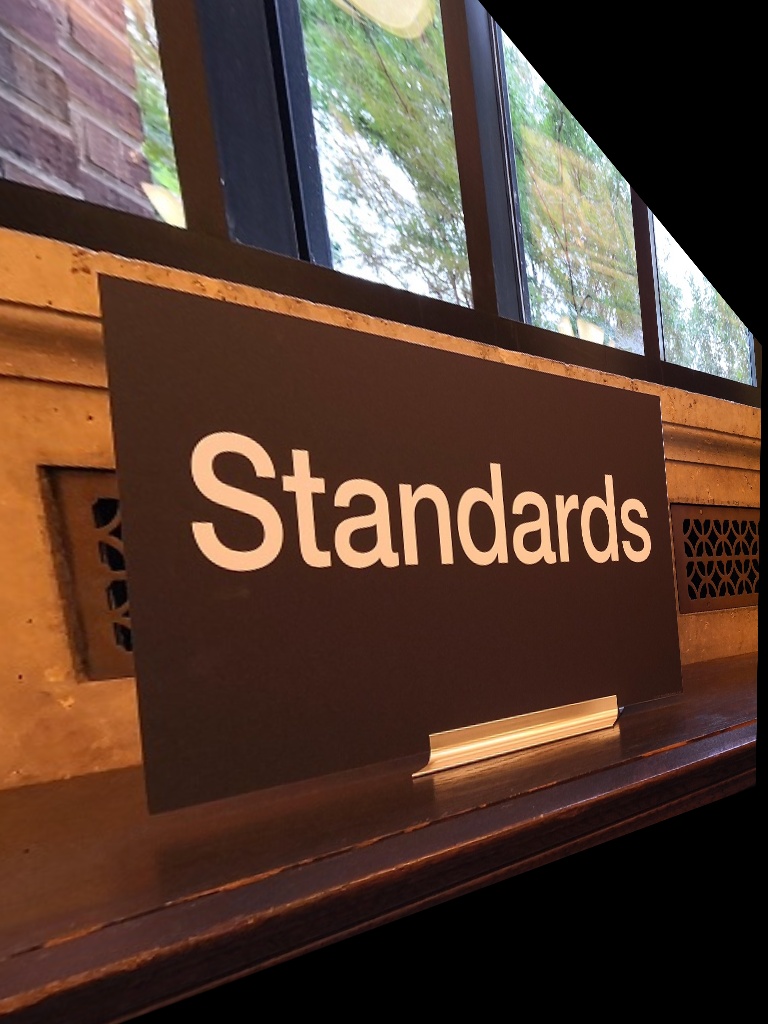}&\includegraphics[width=0.18\textwidth,clip]{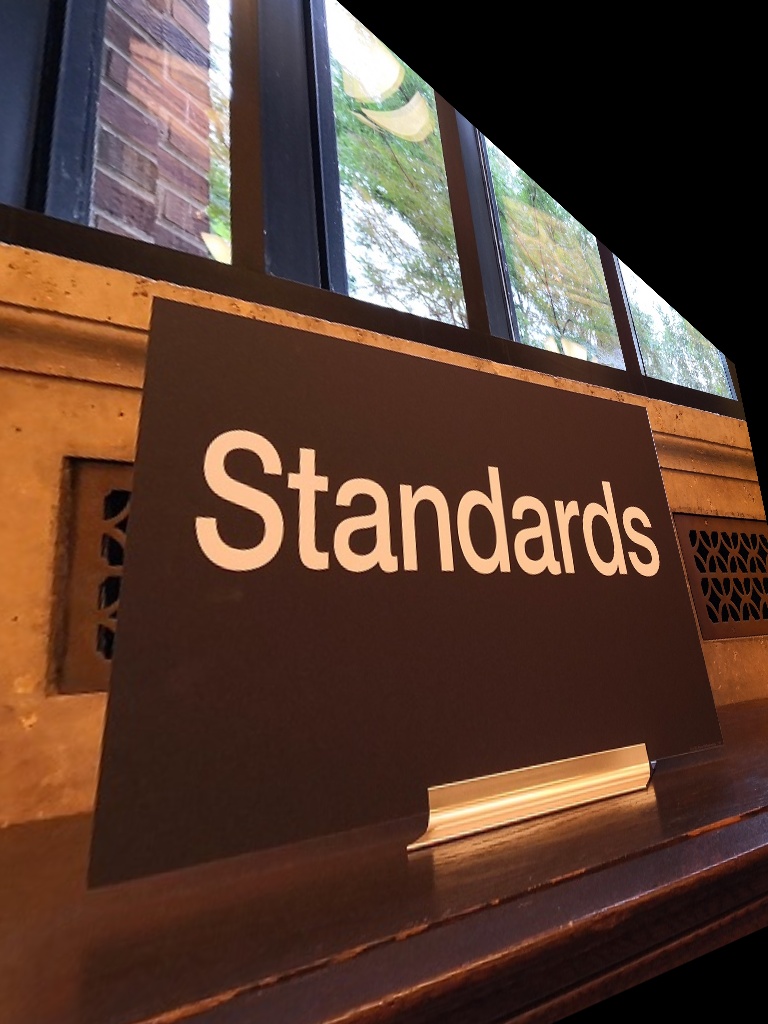}&\includegraphics[width=0.18\textwidth,clip]{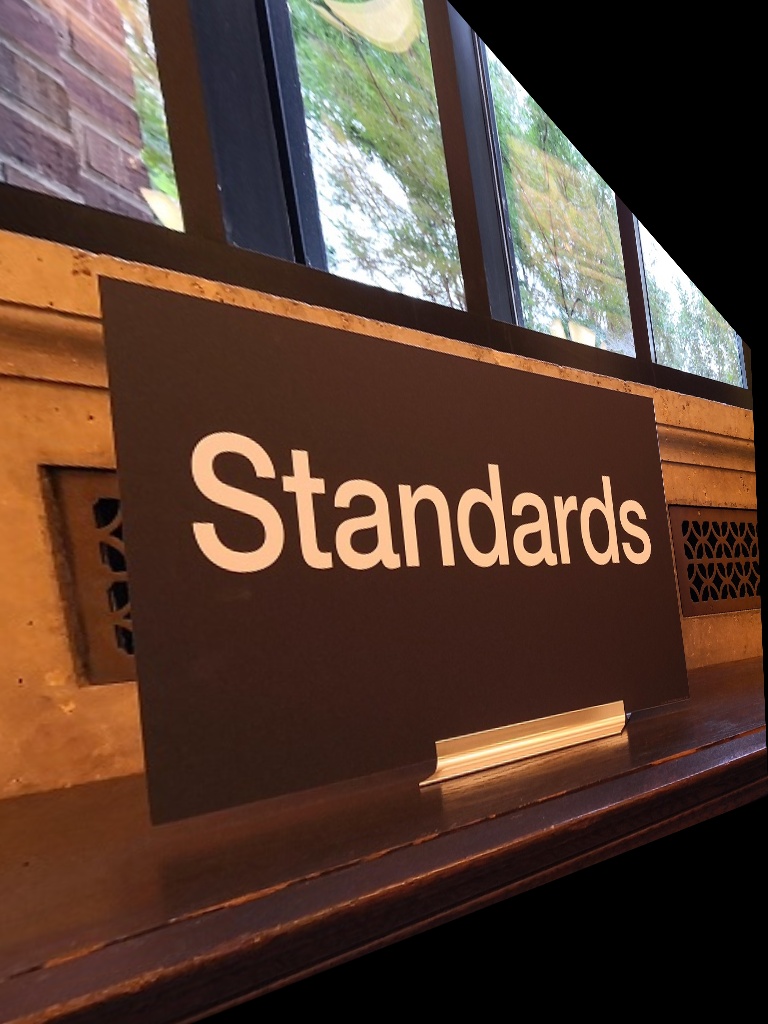}&\includegraphics[width=0.18\textwidth,clip]{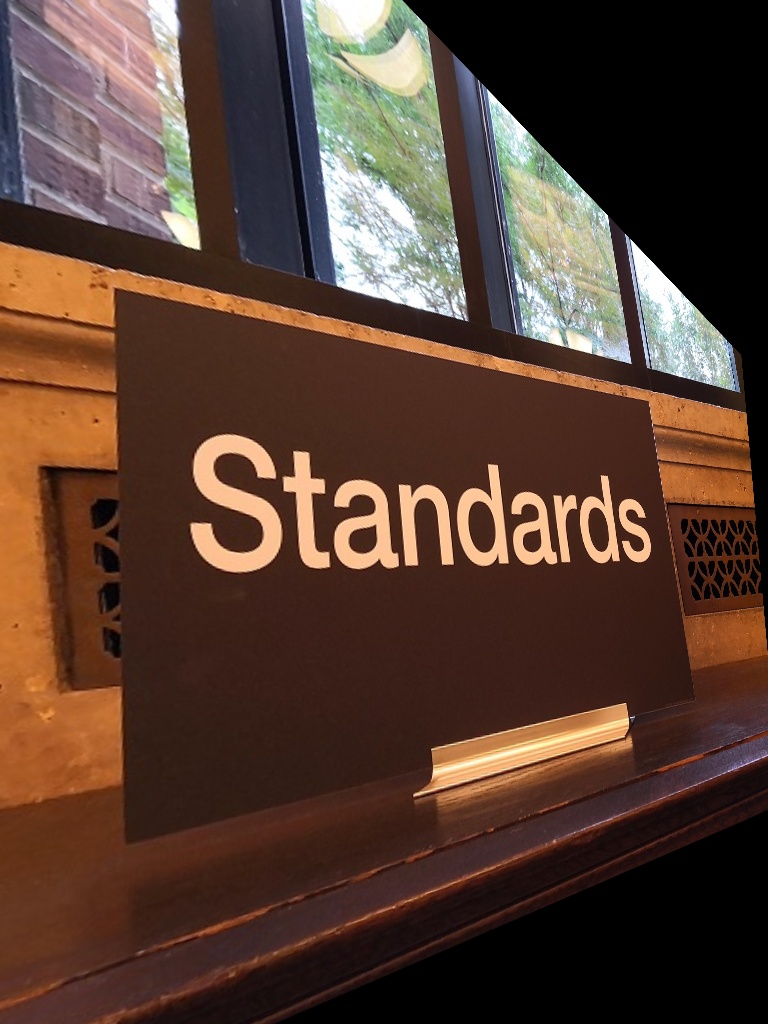}\\
 \includegraphics[width=0.18\textwidth,clip]{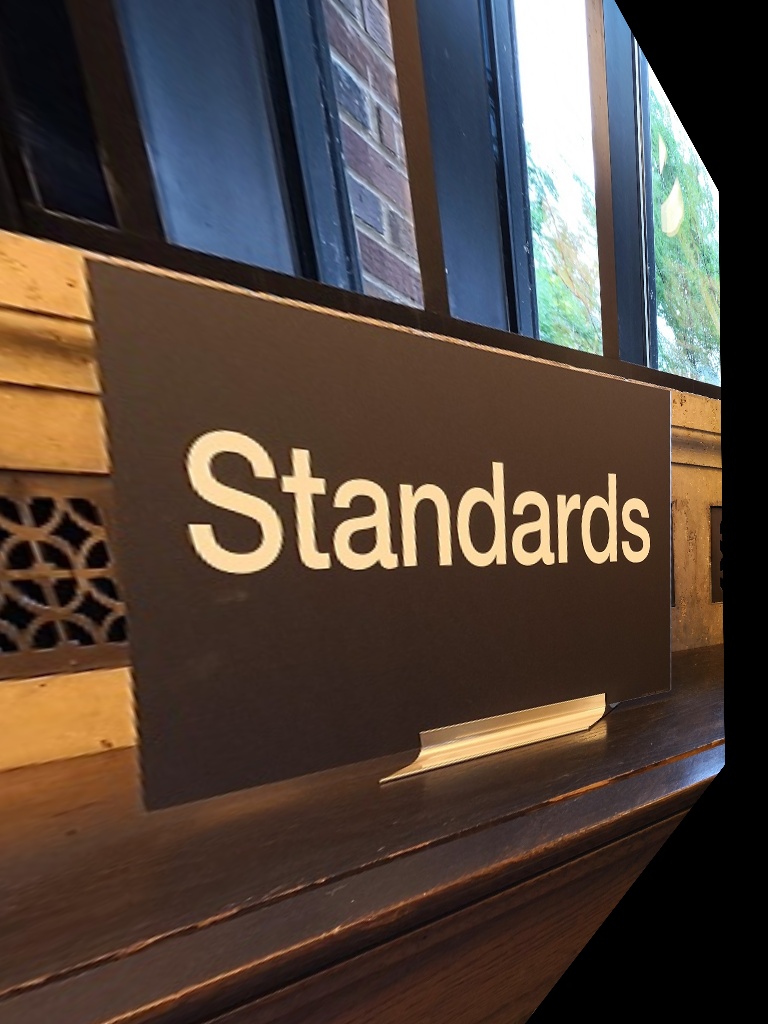}&\includegraphics[width=0.18\textwidth,clip]{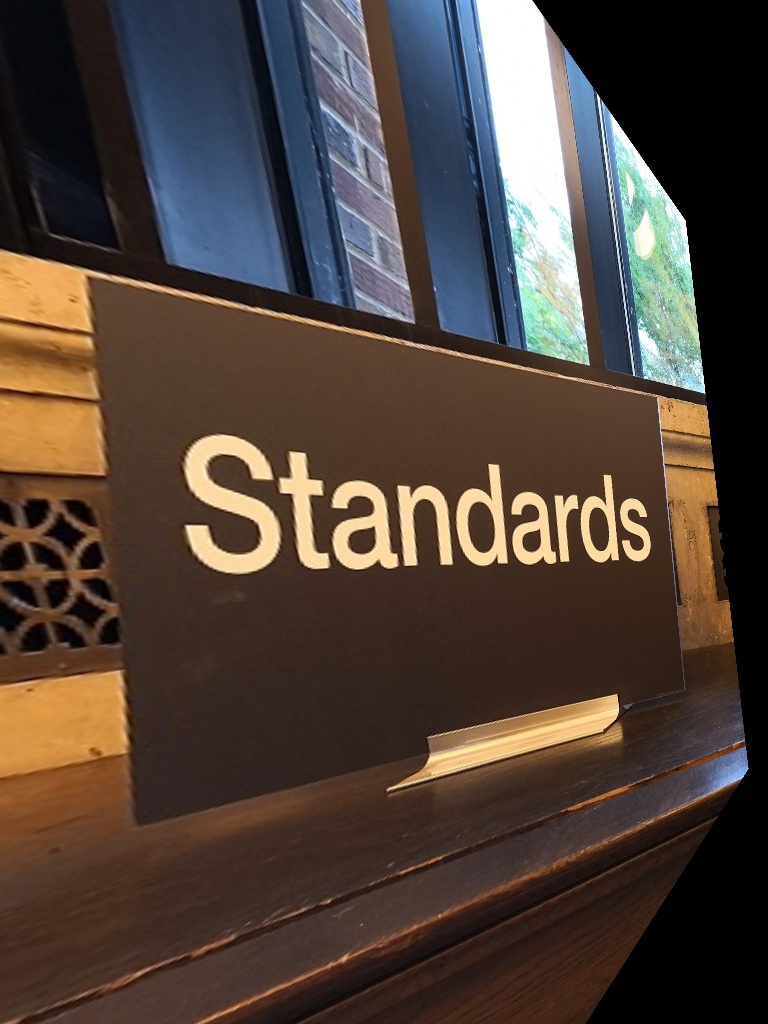}&\includegraphics[width=0.18\textwidth,clip]{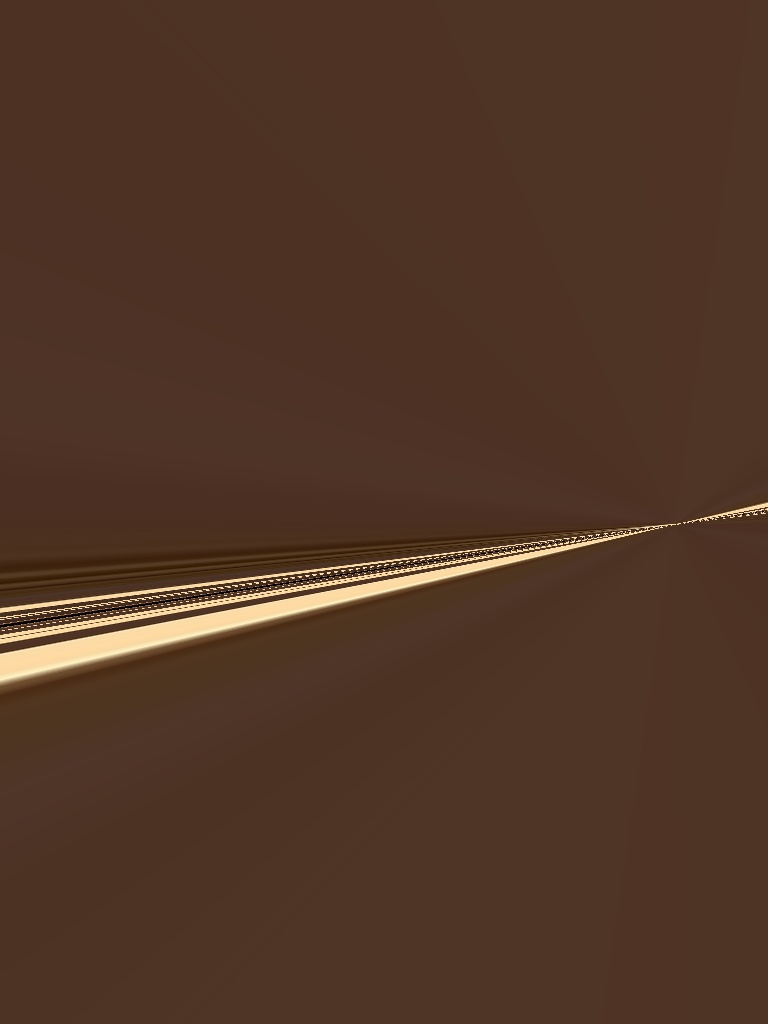}&\includegraphics[width=0.18\textwidth,clip]{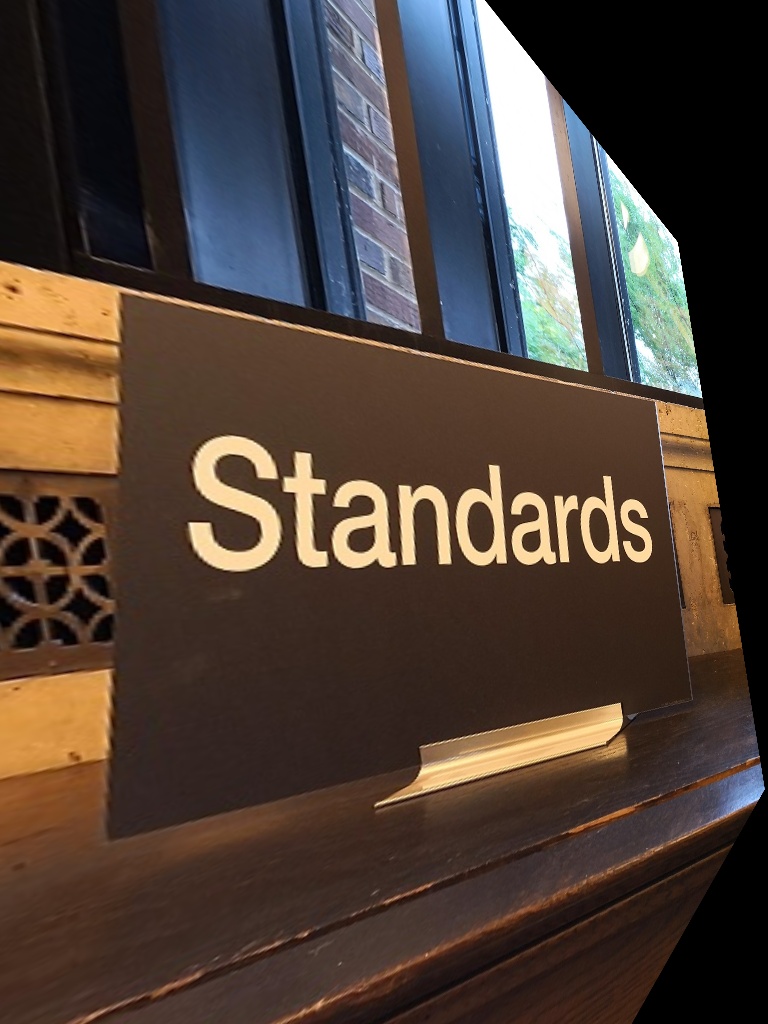}\\
 \includegraphics[width=0.18\textwidth,clip]{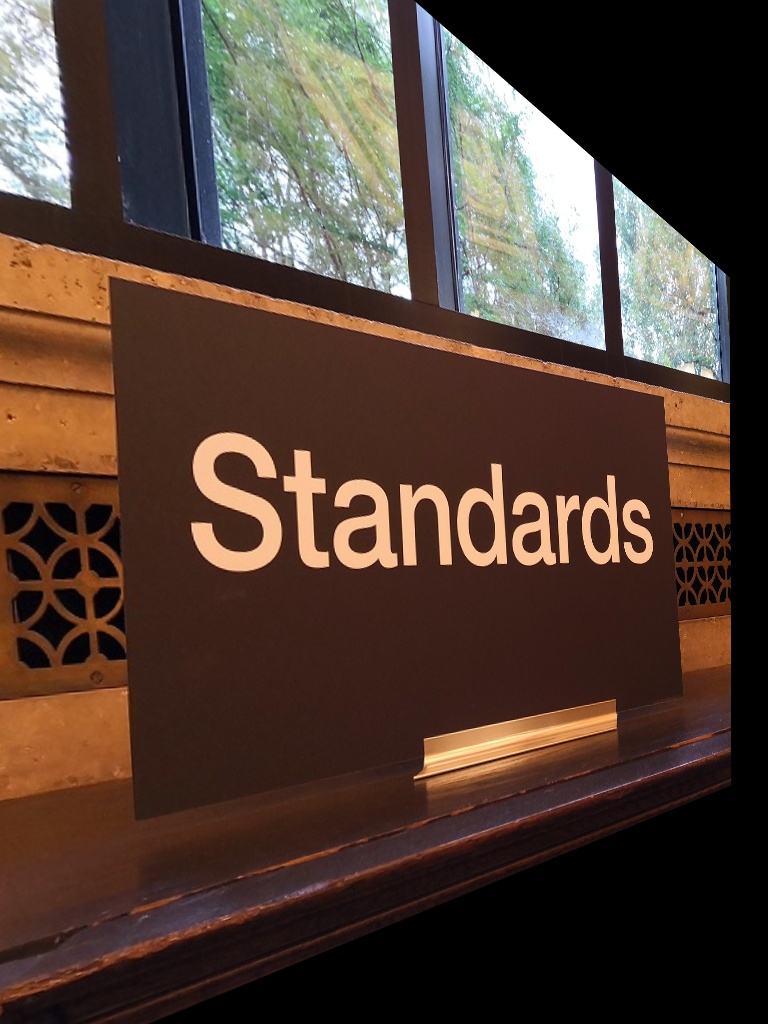}&\includegraphics[width=0.18\textwidth,clip]{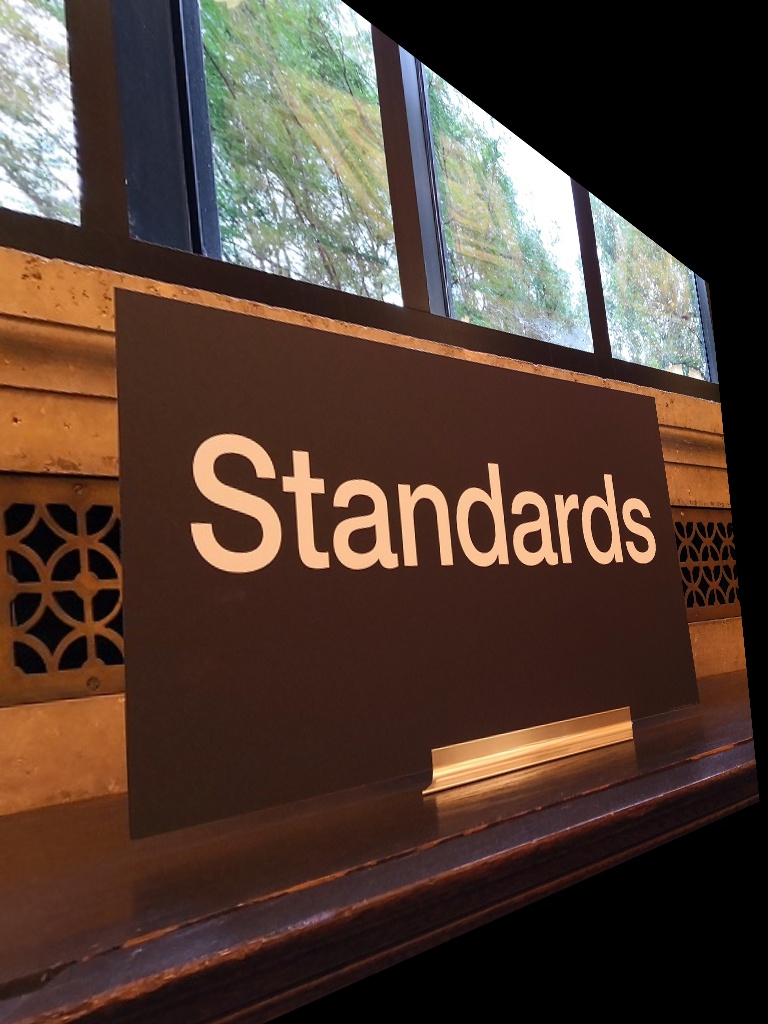}&\includegraphics[width=0.18\textwidth,clip]{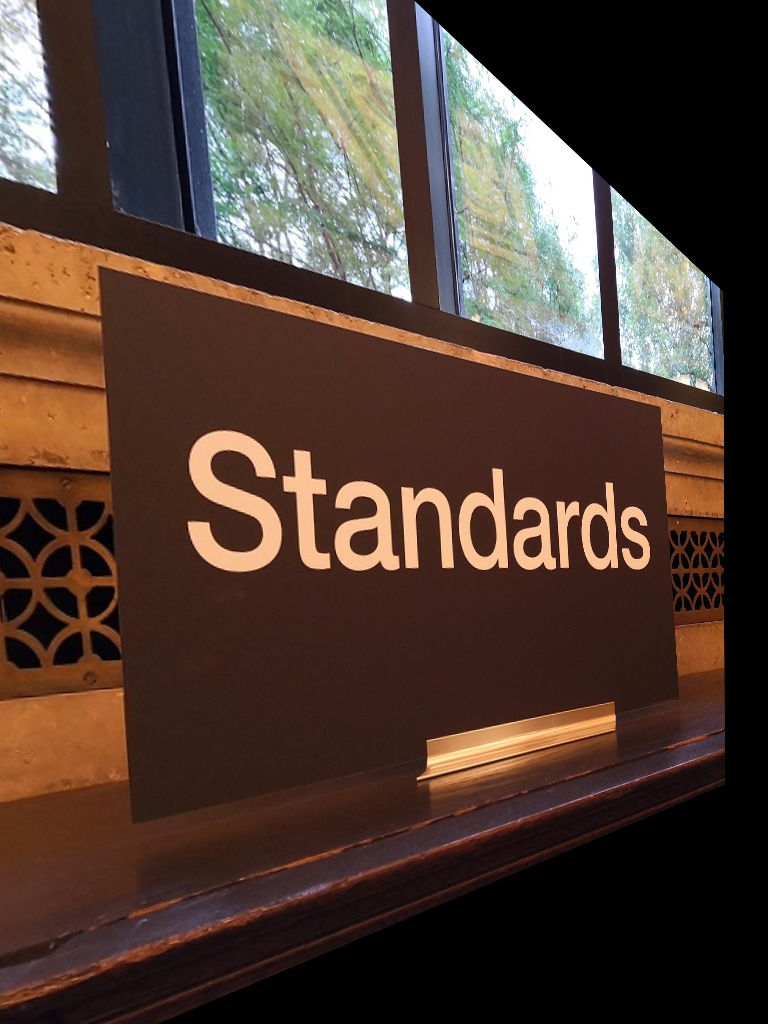}&\includegraphics[width=0.18\textwidth,clip]{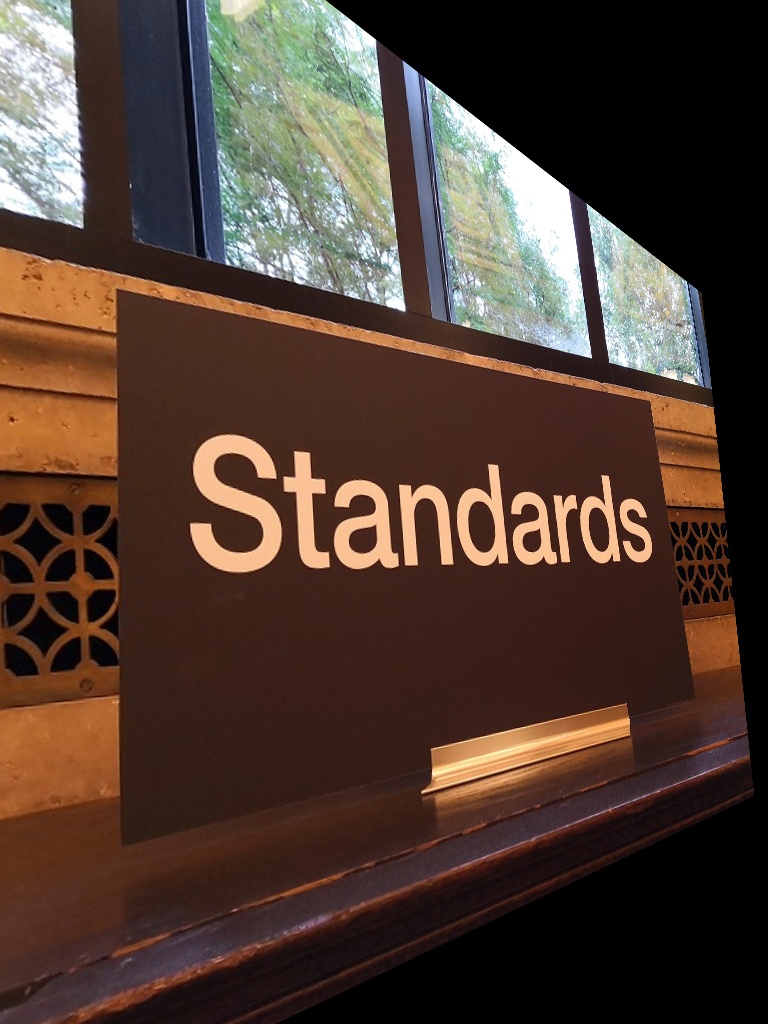}
 \end{tabular}
 \caption{The homography transformation results by using the corresponding points obtained in Figs.~\ref{compareFig} and~\ref{siftcps}.} \label{sifthm}
 \end{figure}

 Further, by using equation~(\ref{arc-p}), we can also obtain the centro-affine arc-lengths for the boundaries of the letters in the label ``Standards'', which are shown as in Table~\ref{length}.
 The centro-affine arc-lengths for the letters in the $i$th image of Fig.~\ref{orspline1-4} are listed in the $i$th row of Table~\ref{length}.
 It is easy to see that there is a strong correspondence between the centro-affine arc-lengths of the corresponding objects.
 Furthermore, we calculate the correlation coefficients between them, as shown in Table~\ref{cor}, which indicate that they are almost same.

In projective geometry, a homography is an isomorphism of projective spaces,
 induced by an isomorphism of the vector spaces from which the projective spaces derive.
 Estimating the~2D homography (or projective transformation) from a pair of images is a fundamental task in computer vision.
 Now let us find the homography transformations between these images and compare the results with the corresponding points obtained by the SURF,
 SIFT, ASIFT, and our centro-affine method.
 Those corresponding points data sets in Figs.~\ref{compareFig} and~\ref{siftcps} are used to transform the other three images to the second one in Fig.~\ref{or1-4},
 the final results are shown in Fig.~\ref{sifthm}.
 The transformation results derived by the corresponding points of the centro-affine invariant method are shown at in the first column; the second column shows the results by the SURF method; the third column shows the results by the SIFT method; the last column shows the results by the ASIFT method. (The second figure in the SIFT column is blank due to incorrect matching points; see the second image of the SIFT row in Fig.~\ref{siftcps}.)

\begin{rem}When the images have less texture complexity and color diversity, it is more challenging to extract and describe the feature points through the SIFT, SURF, and ASIFT methods. The centro-affine invariant method can match the corresponding points by the closed boundary, and hence is less affected by texture and color. In view of Figs.~\ref{compareFig}, \ref{siftcps} and \ref{sifthm}, we find that between the images under the large scale transformation of the camera, the centro-affine invariant method for finding the corresponding points, using the corresponding boundaries, offers certain advantages over SURF, SIFT and ASIFT. However, if the object itself admits affine symmetries, for example, the ``o'' inside the letter ``d'', it is not so easy to find the proper corresponding points; this defect can be observed in Fig.~\ref{compareFig}.
 \end{rem}

\section{Concluding remarks}\label{section6}

 In this paper, we have investigated the heat flow in centro-affine geometry and applications of differential centro-affine invariants in edge matching. Differential invariants and their algebraic relations for curves in centro-affine geometry are easily obtained by applying the equivariant moving frame method~\cite{fo-2}.
 A classification of curves with constant centro-affine curvature has been provided.
 More interestingly, we have shown that the heat flow in centro-affine geometry is equivalent to the first-order inviscid Burgers' equation, in contrast to the equations governed by the heat flows in Euclidean, affine, and centro-equi-affine geometries,
 which are nonlinear second-order parabolic equations.
 Thus, the evolution process for the centro-affine heat flow is described by solving the inviscid Burgers equation through the method of characteristics. An interesting question, to be explored later, is what the presence of shock waves in the solution might mean for the corresponding curve evolution.
 In addition, an application of centro-affine invariants to edge matching is presented. It turns out that the resulting method offers
 certain advantages over other well-used methods.

 In conclusion, we would like to mention further possible issues relating to this work.
\begin{itemize}\itemsep=0pt
\item The Gaussian kernel, while being one of the most used in image analysis,
 has several undesirable properties, principally when applied to planar curves.
 One of these is that the filter is not intrinsic to the curve.
 This can be remedied by replacing the linear heat equation by a geometric heat equation.
 In particular, if the Euclidean geometric heat flow is used, a scale space invariant to rotations and translations is obtained, while the \mbox{(equi-)}affine version leads to an(equi-)affine invariant multi-scale representation; see~\mbox{\cite{ost-1,tot}} for general results and recent developments. It would thus be of interest to apply the centro-affine heat flow analyzed here to construct a corresponding scale space.
\item As shown here, the differential invariants for curves in centro-affine geometry can be used to study edge matching in images.
 A natural question arises: whether we can use differential invariants for surfaces in centro-affine geometry to study the edge matching of three-dimensional images.
 In fact, the differential invariants for the equi-affine group acting on image volumes have been constructed in~\cite{tot}.
\end{itemize}

\subsection*{Acknowledgements}
The authors thank the referees for their valuable comments and suggestions. The work of C.Z.~Qu is supported by the NSF-China grant-11631007 and grant-11971251.

\pdfbookmark[1]{References}{ref}
\LastPageEnding

\end{document}